\documentclass[12pt]{article}

\usepackage[dvips]{graphicx}
\usepackage{makeidx}
\usepackage{amsmath}
\usepackage{amsfonts}
\usepackage{amssymb}
\usepackage{xcolor}

\usepackage[utf8]{inputenc}
\usepackage[T1]{fontenc}
\usepackage{lmodern}

\usepackage{color}

\newtheorem{lem}{Lemma}

\newtheorem{assume}{Assumption}
\newtheorem{cor}{Corollary}
\newtheorem{rmk}{Remark}
\newtheorem{dfn}{Definition}
\newtheorem{prop}{Proposition}
\newtheorem{example}{\it Example}

\newfont{\pseudocode}{cmtt10}

\newcommand{\dist}{\mathrm{dist}}

\newcommand{\pd}{{\partial}}

\newcommand{\Real}{\mathbb{R}}

\newcommand{\Numbers}{\mathbb{Z}}

\renewcommand{\Real}{\mathbb{R}}

\begin{document}

\title{Lyapunov-like Conditions of Forward Invariance and Boundedness for a Class of Unstable Systems
}
\author{A. Gorban\thanks{Dept. of Mathematics, University of Leicester, Leicester, LE1 7RH, UK (ag153@le.ac.uk)}, I. Tyukin\thanks{{\bf Corresponding author}. Dept. of Mathematics, University of Leicester, Leicester, LE1 7RH, UK (I.Tyukin@le.ac.uk);}, E. Steur\thanks{Laboratory for Perceptual Dynamics, KU Leuven, Tiensestraat 102, 3000 Leuven, Belgium (erik.steur@ppw.kuleuven.be,  cees.vanleeuwen@ppw.kuleuven.be)}, and H. Nijmeijer\thanks{Dept. of Mechanical Engineering, Eindhoven University of Technology, P.O. Box 513 5600 MB,  Eindhoven, The Netherlands, (e-mail: h.nijmeijer@tue.nl)}}
\date{}
\maketitle

\begin{abstract}
We provide Lyapunov-like characterizations of boundedness and
convergence of non-trivial solutions for a class of systems
with unstable invariant sets. Examples of systems to
which the results may apply include interconnections of stable
subsystems with one-dimensional unstable dynamics or critically stable
dynamics. Systems of this type arise in  problems of nonlinear
output regulation, parameter estimation and adaptive control.
In addition to providing boundedness and convergence criteria
the results allow to derive domains of initial conditions
corresponding to solutions leaving a given neighborhood of the
origin at least once. In contrast to other works addressing
convergence issues in unstable systems, our results require
neither  input-output characterizations for
the stable part nor estimates of convergence rates.
The results are illustrated with examples, including the analysis of phase synchronization of neural oscillators with heterogenous coupling.
\end{abstract}


\section{Introduction}\label{sec:intro}

 Methods and tools for the analysis of asymptotic properties of solutions of ordinary differential equations are
important components of modern control theory. Even though
the problems of control are often viewed as that of synthesis
rather than analysis, the latter crucially affects the former.
Indeed, in order to be able to specify feasible goals of
synthesis, e.g. forward-completeness, state boundedness,
asymptotic convergence of solutions to a region in the state
space etc., one needs to understand how these properties depend
on the system parameters and controls.

The majority of the analysis techniques in control, and hence
methods for systems design,  rely upon the assumption that desired
motions in the system are stable in the sense of Lyapunov
\cite{Lyapunov:1892}. Let us briefly recall this and other
related notions from the domain of dynamical systems,  and also introduce notational conventions used throughout the manuscript.

\subsection{Notation and basic notions}

The following notational conventions are used throughout the
paper.   Let $\mathcal{D}$ be an open set in $\Real^n$.
 The symbol $\mathcal{C}^k(\mathcal{D})$ denotes the space of
functions that are  at least $k$ times differentiable in
$\mathcal{D}$;  $\overline{\mathcal{D}}$ denotes the closure of
$\mathcal{D}$; $\|\cdot\|$ stands for the Euclidian norm. Let $\mathcal{S}$ be a subset of $\Real^n$, and $x\in\Real^n$, then $\dist(\mathcal{S},x)=\inf_{x'\in\mathcal{S}}\|x-x'\|$. By
$\mathcal{K}_0$ we denote   the set of all non-decreasing
continuous functions $\kappa: \Real_{\geq 0}\rightarrow
\Real_{\geq 0}$ such that $\kappa(0)=0$; $\mathcal{K} \subset \mathcal{K}_0$
is the subset of strictly increasing functions, and $\mathcal{K}_\infty \subset \mathcal{K}$
consists of functions from $\mathcal{K}$ with infinite limit: $\lim_{s\rightarrow\infty}\kappa(s)=\infty$.
  Consider a non-autonomous system $\dot{x}=f(x,p,t)$,  where $f:\Real^n\times\Real^d\times\Real\rightarrow\Real^n$
is continuous, $p\in\Real^d$ is the vector of parameters, and $f(\cdot,p,t)$ is locally Lipschitz; $x( \cdot \ ;t_0,x_0|p)$ stands for
the unique maximal solution of the initial value problem: $x(t_0;t_0,x_0|p)=x_0$. In cases when no confusion arises, we will refer to these solutions as $x(\cdot;t_0,x_0)$,  $x(\cdot;x_0)$,  or simply $x(\cdot)$. Solutions of the initial value problem above at $t$ are denoted as
$x(t;t_0,x_0)$,  $x(t;x_0)$,  or $x(t)$ respectively.
We always separate by the semicolon the symbol of the independent (time)
variable from symbols of other variables (initial data or parameters).

Let us start with the classical notion of invariance of a set. Let $\mathcal{D}$ be an open subset of $\Real^n$, and consider systems represented by differential equations
$\dot{x}=f(x)$    in the domain $\mathcal{D}$. The right-hand side, $f(x)$, is
assumed to be a  {\em locally Lipschitz} vector-field on $\mathcal{D}$. In this case, for any initial condition $x(0)=x_0, x_0 \in \mathcal{D}$
the system has a solution   $x(\cdot;x_0)$   defined on
 a time interval $(-\tau,\tau)$ where $\tau>0$ may depend on $x_0$.
A set $\mathcal{S}$, {  $\mathcal{S}\subset\mathcal{D}$,} is {\em {  forward} invariant} (w.r.t. the system dynamics)
if for every $x_0 \in \mathcal{S}$,    $x(\cdot;x_0)$  is defined on
$[0,\infty)$ and  $x(t;x_0)\in \mathcal{S}$ for all $t>0$.
$\mathcal{S}$  is {\em invariant} if for every $x_0 \in \mathcal{S}$ the solution   $x(\cdot;x_0)$   is defined on $(-\infty,\infty)$, and
$x(t;x_0)\in \mathcal{S}$ for all $t\in\Real$. Unions and intersections of a family of (forward) invariant
sets are (forward) invariant.


A closed invariant set $\mathcal{S}\subset\mathcal{D}$ is a {\it weakly attracting} set if there exists
a set $\mathcal{V}\subset \mathcal{D}$ with strictly
positive measure such that for all $x_0\in\mathcal{V}$ the solution $x(\cdot;x_0)$   is defined on $[0,\infty)$ and
the following holds: $\lim_{t\rightarrow\infty}\mathrm{dist}(\mathcal{S},x(t;x_0))=0$
\cite{Milnor_1985}. The  set $\mathcal{V}$ is not necessarily a neighborhood of
$\mathcal{S}$. The set is {\it attracting} if
$\mathcal{V}$ is a neighborhood of $\mathcal{S}$, and
$\mathcal{V}$ is forward invariant. A closed invariant set
$\mathcal{S}\subset \mathcal{D}$ is {\it stable in the sense of
Lyapunov} if for any neighborhood $\mathcal{V}$ of
$\mathcal{S}$ there exists a {  forward} invariant neighborhood
$\mathcal{W}\subset\mathcal{D}$ of $\mathcal{S}$ such that
$\mathcal{W}\subset\mathcal{V}$ \cite{Zubov:1964}. In other
words, a set that is stable in the sense of Lyapunov has a {\it
fundamental base} of {  forward} invariant neighborhoods. (A collection
$\mathcal{U}_\mathcal{S}$ of all neighborhoods of $\mathcal{S}$
is called a {\it neighborhood system} of $\mathcal{S}$. A
subcollection
$\mathcal{B}_\mathcal{S}\subset\mathcal{U}_\mathcal{S}$ is a
{\it fundamental base} of system $\mathcal{U}_\mathcal{S}$ iff
every element of $\mathcal{U}_\mathcal{S}$ contains at
least one element of $\mathcal{B}_\mathcal{S}$.)

For non-compact sets $\mathcal{S}\subset\mathcal{D}$ it may be useful to distinguish the notion of
Lyapunov stability from the notion of {\em uniform Lyapunov stability} that is defined with uniform neighbourhoods \cite{England:1967}.

Various extensions of stability of sets are proposed for nonautonomous systems too \cite{Kloeden2011}.
For these systems we need the notion of forward invariance of sets in the state space. Consider systems  $\dot{x}=f(x,t)$ in a domain $\mathcal{D}\times\Real \subset \Real^{n+1}$, where the vector-field $f:\mathcal{D}\times\Real \rightarrow \Real^{n}$ is continuous, and $f(\cdot,t)$ is {\em locally Lipschitz} uniformly in $t$.  A set $\mathcal{S} \subset \mathcal{D}$ is { $t_0$-forward} invariant w.r.t. dynamics if for given $t_0\in \Real$  and every
$x_0 \in \mathcal{S}$ the solution   $x(\cdot;t_0, x_0)$ is defined on $[t_0,\infty)$ and $x(t;t_0, x_0) \in \mathcal{S}$
for $t \geq t_0$. If it is $t_0$-forward invariant for all $t_0$ then we call it forward invariant. In this work, we use systems of nested forward invariant sets to characterize
the attractivity of solutions that are not stable in the classical senses.

The notion of Lyapunov stability and analysis methods that are
based on this notion are proven successful in a wide range of
engineering applications (see e.g. \cite{Nijmeijer_90},
\cite{Isidory}, \cite{Ljung_99}, \cite{Narendra89} is a
non-exhaustive list of references).
The popularity and success of the concept of Lyapunov stability
resides, to a substantial degree, in the convenience and
utility of the method of Lyapunov functions for
assessing asymptotic properties of solutions of ordinary
differential equations. Instead of deriving the solutions
explicitly it suffices to solve an algebraic inequality
involving partial derivatives of a given Lyapunov candidate
function.

As the methods of control expand from purely engineering
applications into a wider area of science, there is a need for
maintaining behavior that fails to obey the usual requirement of
Lyapunov stability. There are numerous examples of systems
possessing Lyapunov-unstable, yet attracting, invariant sets
\cite{Andronov:1973},  e.g., in the domains of aircraft dynamics and design
of synchronous generators \cite{Bautin:1990} (pp. 313--356).
Even though solutions of these systems may not always be
uniformly asymptotically stable, they are required to be
bounded and converging  to some specified areas in the system
state space. Finding rigorous and tight criteria for asymptotic
convergence to Lyapunov-unstable invariant sets, however, is a
non-trivial problem.

\subsection{Motivating Examples}

Let us start from examples illustrating non-trivialities of
asymptotic behavior of solutions in systems with unstable
invariant sets.
\begin{example}\label{example:systems_toy}\normalfont Consider the following systems:\\
\begin{subequations}
\begin{minipage}{0.3\textwidth}
\begin{equation}
\left\{ \begin{array}{ll}\dot{x}&=-x + \lambda\\ \dot{\lambda}&=-\gamma |x|^3\end{array} \right.\label{eq:simple_system:1}
\end{equation}
\end{minipage}
\begin{minipage}{0.3\textwidth}
\begin{equation}
\left\{ \begin{array}{ll}\dot{x}&=-x^2 + \lambda\\ \dot{\lambda}&=-\gamma |x|^3\end{array} \right.\label{eq:simple_system:2}
\end{equation}
\end{minipage}
\begin{minipage}{0.3\textwidth}
\begin{equation}
\left\{ \begin{array}{ll}\dot{x}&=-x^3 + \lambda\\ \dot{\lambda}&=-\gamma |x|^3.\end{array} \right.\label{eq:simple_system:3}
\end{equation}
\end{minipage}
\end{subequations}\\

\noindent Phase plots of these systems are shown
in Fig. \ref{fig:example:attractors}.
\begin{figure}[th!]
\centering
\begin{minipage}[h]{0.3\linewidth}
\centering
\includegraphics[width=\textwidth]{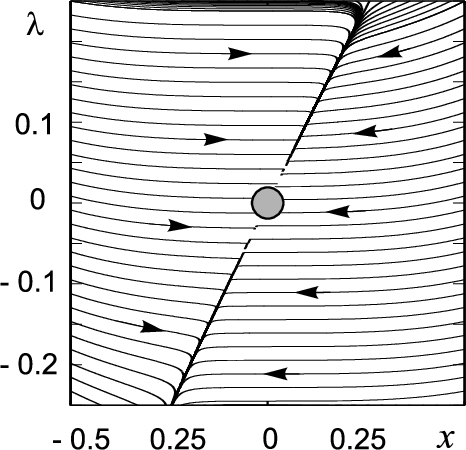}
\vspace{3mm}
a
\end{minipage}
\hspace{3mm}
\begin{minipage}[h]{0.3\linewidth}
\centering
\includegraphics[width=\textwidth]{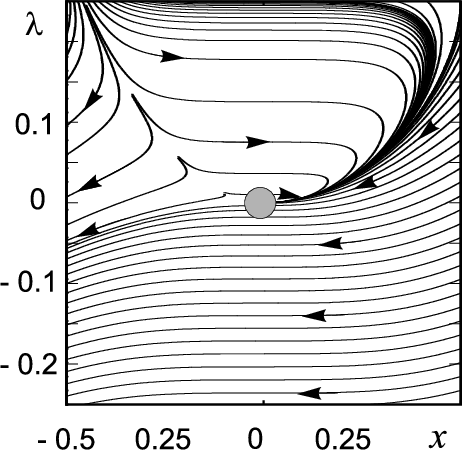}
\vspace{3mm}
b
\end{minipage}
\hspace{3mm}
\begin{minipage}[h]{0.3\linewidth}
\centering
\includegraphics[width=\textwidth]{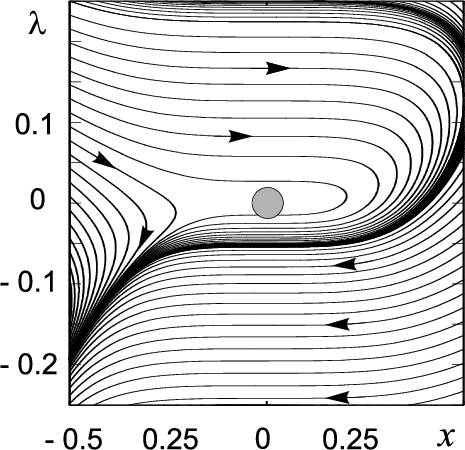}
\vspace{3mm}
c
\end{minipage}
\vspace{-3mm}
\caption{Phase plots of  (\ref{eq:simple_system:1}) (panel $a$), (\ref{eq:simple_system:2}) (panel  $b$), and (\ref{eq:simple_system:3}) (panel  $c$) for $\gamma=0.5$. }\label{fig:example:attractors}
\end{figure}
Systems (\ref{eq:simple_system:1}), (\ref{eq:simple_system:3})
share apparent similarity in their description. Indeed, $x=0$
is the unique asymptotically stable equilibrium of their first
equations at $\lambda=0$, and their second equations governing
the dynamics of variable $\lambda$ are identical. System
(\ref{eq:simple_system:2}), to the contrary, is fundamentally
different from (\ref{eq:simple_system:1}),
(\ref{eq:simple_system:3}): the equilibrium $x=0$ of its first
equation is not stable at $\lambda=0$.  Thus one could expect
qualitative similarity of the dynamics of
(\ref{eq:simple_system:1}) and (\ref{eq:simple_system:3}),
whereas no such similarity is expected between
(\ref{eq:simple_system:2}) and (\ref{eq:simple_system:1}).
However, as the phase plots suggest (see Fig.
\ref{fig:example:attractors}), the origins of
(\ref{eq:simple_system:1}), (\ref{eq:simple_system:2}) are weak
attractors, and the origin of (\ref{eq:simple_system:3}) is not
an attractor. This can be demonstrated, for example, by using the
following singular transformation: $(x,\lambda)\mapsto
(\rho,\varphi)$,  $x=\rho \cos{\varphi}$,
$\lambda=\rho^2\sin(\varphi)$ (for details see \cite{Bogoyavlenskii:1980}).
Therefore, our initial naive intuition about qualitative properties of solutions of
(\ref{eq:simple_system:1})--(\ref{eq:simple_system:3}) failed.
This motivates the necessity for having rigorous, simple and
efficient criteria for assessing asymptotic properties of
solutions in systems with unstable invariant sets.
\end{example}

In pursuing this goal we do not wish to attempt addressing the
issue in full generality, viz. for the widest class of systems
possible. Instead we focus on a particular family of equations
that occur naturally in a range of control, optimization,
estimation and modelling problems. Examples of systems from
this family are provided below.

\begin{example}\label{example:0}\normalfont{\it A  Network of phase oscillators} provides an example of systems in which unstable attractors are shown (numerically) to exist and sometimes prevail \cite{PRL:2002:Timme}, \cite{Chaos:2003:Timme}. One of the simplest instances of such networks is a network of three phase oscillators of which the phase differences are governed by the following set of equations (see Section \ref{subsec:example:phase_synch}):
\begin{subequations}
\begin{eqnarray}
\dot{\vartheta}_1&=&  \varepsilon/\pi\sin(\vartheta_2)^2\label{eq:cell:2:p}\\
\dot{\vartheta}_2&=& - \varepsilon_1/2\sin(2 \vartheta_2)+\varepsilon/\pi(\sin(\vartheta_1)^2+\sin(\vartheta_2)^2)\label{eq:cell:3:p}.
\end{eqnarray}
\end{subequations}
Variables  $\vartheta_1$,   $\vartheta_2$ denote the phase
differences, and $\varepsilon$, $\varepsilon_1$ are parameters
representing the coupling strengths between oscillators. It is
clear that the origin of (\ref{eq:cell:2:p}),
(\ref{eq:cell:3:p}) is unstable in the sense of Lyapunov. This rules out explicit application
of Lyapunov functions-based approaches for deriving conditions
of synchronization in such systems. The questions nevertheless
are: for which values of $\varepsilon,\varepsilon_1$ phase
synchronization will occur,  how large is the domain of initial
conditions leading to such synchronous state, and how  does it
depend on $\varepsilon,\varepsilon_1$?
\end{example}

\begin{example}\label{example:1}\normalfont  {\it A universal adaptive stabilizer of nonlinear systems}  (in presence
of uncertainties). The studied question is as follows. Consider the following system
\begin{equation}\label{eq:problem:1}
\begin{split}
\dot{x}&=f(x,t)+g(x,t)u(t), \ f:\Real^n\times\Real\rightarrow\Real^n, \ g:\Real^n\times\Real\rightarrow\Real^{n\times m}\\
y&=h(x), \ h:\Real^n\rightarrow\Real^d,
\end{split}
\end{equation}
where $y$ is the output, $u:\Real\rightarrow\Real^m$ is a control input which may depend on the current value of $y(t)=h(x(t))$ and  on time $t$ directly. The functions $f(\cdot,\cdot),g(\cdot,\cdot)$ are not known explicitly.
What a priori information about the system should be made available in order to derive a control input $u(\cdot)$ stabilizing the zero solutions of (\ref{eq:problem:1})? The question has been answered in \cite{Dyn_Con:Ilchman:97}, \cite{CDC:Pomet:1992} (see also \cite{Automatica:Algower:1997,IEEE_TAC:Ilchmann:1999,IEEE_TAC:Ilchmann:2003} for related work). Assume that there exists a matrix $K\in\Real^{m \times d}$ such that the zero solution of
$\dot{x}=f(x,t)+g(x,t)K h(x)$ is exponentially stable. The following system  with the stabilizer was constructed:
\begin{equation}\label{eq:problem:1:2}
\begin{split}
\dot{x}&=f(x,t)+g(x,t)\beta(\gamma(\lambda))h(x)\\
\dot{\lambda}&=\|h(x)\|^p,
\end{split}
\end{equation}
where $\beta(\cdot):\Real\rightarrow\Real^{m\times d}$ is a special function with dense image in $\Real^{m\times d}$ and $\gamma:\Real\rightarrow\Real$ is a special monotone function of which the growth rate decays to zero as $\lambda\rightarrow\infty$. It was proven that the $x$-component of solutions, $x(\cdot)$, converges to
the origin as $t\rightarrow\infty$. (For the details of $\beta (\cdot)$ and $\gamma (\cdot)$ construction and specific conditions on functions $f(\cdot, \cdot)$, $g(\cdot, \cdot)$ and $h(\cdot)$ see original papers.) Nevertheless, solutions of the extended system (\ref{eq:problem:1:2}) are not uniformly asymptotically stable (cf. \cite{SICON:Townley}).
\end{example}

Despite the  problems described in Examples \ref{example:0} and
\ref{example:1} arise in different subject areas, they are
inherently similar. In both cases we have to deal with systems
composed of an ``attracting'' subsystem coupled with a
``wandering'' one. The attracting subsystem has an attracting
invariant set in its state space, and solutions of the
wandering subsystem unidirectionally evolve along a certain
path.   The general description  of the composed system is provided below.

\subsection{Class of Systems}

We will focus on systems whose dynamics can be described by the system of ordinary differential equations:
\begin{subequations}
\begin{equation}
\left\{ \begin{array}{ll}\dot{x}&=f(x,\lambda,t)\\
\dot{\lambda}&=g(x,\lambda,t),
\end{array} \right.\label{eq:problem}
\end{equation}
\noindent   where the vector-fields $f: \Real^n\times\Real\times\Real\rightarrow \Real^n$, $g:\Real^n\times\Real\times\Real\rightarrow\Real$ are continuous,  $f(\cdot,\cdot,t)$, $g(\cdot,\cdot,t)$ are locally Lipschitz uniformly in $t$,
and  $g(\cdot,\cdot,\cdot)$ in (\ref{eq:problem}) is of constant sign.
    Eq. (\ref{eq:problem})
describes the coupled system generalizing (\ref{eq:simple_system:1})--(\ref{eq:simple_system:3}), (\ref{eq:cell:2:p}),(\ref{eq:cell:3:p}), and (\ref{eq:problem:1:2}).

Throughout the paper we assume that $(0,0)$ is an equilibrium of (\ref{eq:problem}). Moreover, we assume that the origin $x=0$ is a weak attactor of the $x$-subsystem of (\ref{eq:problem}) for frozen values of $\lambda$ at $\lambda=0$, i.e. the origin $x=0$ of the system
\begin{equation}
\left\{ \begin{array}{ll}\dot{x}&=f(x,\lambda,t)\\  \lambda&={\mathrm{const}}, \lambda \in \Real. \end{array} \right.\label{eq:problemFroz}
\end{equation}
\end{subequations}
is a weak attractor at $\lambda =0$.
{  We also assume that for system  (\ref{eq:problemFroz}) there is a $p>0$ and a set
$\underline{\omega}(p)$ which is forward invariant for all $\lambda\in[0,p]$: if $x_0\in\underline{\omega}(p)$ and  $\lambda\in[0,p]$ then
\begin{equation}\label{eq:pre-bag}
  x(t;t_0,x_0| \lambda)\in\ \underline{\omega}(p)   \ \mbox{for all} \ t_0\in\Real, t\geq t_0.
\end{equation}}
{ In principle, { $\underline{\omega}(p)$} is allowed to coincide with
$\Real^n$. Notice also that since we don't wish to impose any additional
specific constraints (such
as e.g. minimality) the set { $\underline{\omega}(p)$} is
not uniquely determined by the system itself. {  For example, for the system (\ref{eq:problemFroz}) induced by the first equation of (\ref{eq:simple_system:1}) the sets $\underline{\omega}(p)$ can be chosen as intervals $[b_1,b_2]$, $b_1\leq 0$, $b_2\geq p$, as well as $(-\infty,b_2]$, $[b_1,\infty)$, or $(-\infty,\infty)$. For the system  (\ref{eq:problemFroz}) corresponding to (\ref{eq:simple_system:2}) the sets $\underline{\omega}(p)$ are $[0,b_1]$, $b_1\geq \sqrt{p}$.}  For the sake of simplicity
one may ignore references to { $\underline{\omega}(p)$} in the statements of the results,
and assume that { $\underline{\omega}(p)$} coincides with $\Real^n$. On the other hand, as we shall see later,
introduction of { $\underline{\omega}(p)$} enables to produce criteria for checking whether an equilibrium is a
weak attractor or not for dynamical systems (e.g. described by
(\ref{eq:simple_system:2})) in which the zero solution of
(\ref{eq:problemFroz}) at $\lambda=0$ is not stable in the sense of
Lyapunov.}

Further and specific technical assumptions about  $f(\cdot,\cdot,\cdot)$ and $g(\cdot,\cdot,\cdot)$ are proided in Section \ref{sec:Formulation}.

\begin{rmk}\label{rmk:ExtInv} If $\underline{\omega}(p)$ is forward invariant for (\ref{eq:problemFroz}) for all $\lambda\in[0,p]$ then the set $\underline{\omega}(p)$ is also forward invariant w.r.t. equations
$$\dot{x}=f(x,\lambda(t),t)$$
for any piecewise-continuous function  $\lambda(\cdot)$ with values $\lambda(t) \in [0,p]$ (and a discrete set of discontinuity points).
This can be easily proved using approximation of $\lambda(\cdot)$ by piecewise constant functions.
\end{rmk}

In adition to the previous examples, equations (\ref{eq:problem}) describe estimation algorithms in problems of adaptive control and observer design when models of uncertainty are
nonlinearly parameterized, or when the application of standard
techniques is computationally ill-posed
\cite{MMNP:Fairhurst:2010,IFAC_CONGRESS_2011:1}.
They also can be viewed as a prototype for control and
estimation schemes with pre-routing  in the domain of
supervisory control \cite{Morse_95}.

A rather general interpretation of systems (\ref{eq:problem})
is that they govern a class of systems in which inherent
dynamics of an object (first equation and (\ref{eq:problemFroz})) is coupled with the
dynamics of the system's resources (second equation). In this
regards  $g(\cdot,\cdot,\cdot)$ defines the rate of the resource's
consumption, and  $f(\cdot,\cdot,\cdot)$  determines the velocity of the state
$x$ given the available resources $\lambda$ at $t$.

\subsection{Extension of Lyapunov's idea onto unstable sets}

The systems with inherently unstable behavior are important for many applications, including
modelling, control and identification (see, for example, \cite{Rabinovic_2006}, \cite{Rabinovic_2008}, where dissipative
saddles are used to model decision-making sequences, \cite{GOMAN},
where the flutter suppressors were developed,  and \cite{Dyn_Con:Ilchman:97,CDC:Pomet:1992},
where the general problem of universal adaptive stabilization was studied).
Nevertheless, there are limitations restricting further
progress in application of the broader concept of unstable
convergence in these areas. Among these is the lack of a simple analogue of the Lyapunov
method for these, strictly speaking, unstable systems that
would allow to draw conclusions about asymptotic properties of
unstable solutions without the need of solving the equations.
This motivates the focus of our present work.

In this paper we propose an extension of the classical Lyapunov
function method for assessing boundedness and convergence of
motion in dynamical systems with unstable invariant
sets. The class of systems we will consider is given by
(\ref{eq:problem}), and the questions we address below can be
formulated as follows:
\begin{enumerate}
\item Let the origin of (\ref{eq:problem}) be an
    equilibrium. Can we tell (without solving the system)
    if the set is an attractor in some appropriate, e.g.
    Milnor's sense \cite{Milnor_1985}?
\item Pick a point in the system's state space. Is
    it possible to predict (without solving the system) if
    the solution passing through this point is bounded in
    forward time, or does it escape to infinity?
\end{enumerate}
These questions are certainly not original. Algebraic criteria
for checking attractivity of unstable point attractors in a
rather general setting have been proposed in
\cite{Rantzer}, and were further developed in
\cite{Masubuchi:2007,Vadia:2008}. These results apply to
systems in which almost all points in a neighborhood of the
attractor correspond to solutions converging to the attractor
asymptotically. Yet, as can be seen clearly from Fig.
\ref{fig:example:attractors}, this requirement  may not hold
for the class of systems described by (\ref{eq:problem}). On
the other hand techniques which can be used to
address the questions above for equations (\ref{eq:problem}),
such as, e.g.,  \cite{SIAM_non_uniform_attractivity}, lack the convenience of
the method of Lyapunov functions. Further, they require
existence of input-output gains for the stable subsystem.
Hence, developing novel methods to address the issue of convergence to unstable
sets is needed.  These methods, on the one hand, should inherit
the efficiency of Lyapunov analysis in which boundedness of
solutions can be verified by checking a system of
inequalities without involving prior knowledge of the solutions
of the system. On the other hand, these methods should apply to
systems with instabilities such as specified by
(\ref{eq:problem}). In our present contribution we provide a
set of results that can be considered as a possible candidate.

The main idea behind the development of these results can
briefly be summarized as follows. Since we are interested in
the solutions that are not necessarily stable in the sense of
Lyapunov we abandon the concept of neighborhoods from standard
Lyapunov analysis \cite{Lyapunov:1892}, \cite{Zubov:1964}. For
a given invariant set $\mathcal{S}$ of a system, instead of
searching for a fundamental base of {  forward} invariant
neighborhoods of  $\mathcal{S}$, we study existence of
a collection of {   forward} invariant sets
associated with $\mathcal{S}$. These sets are not
necessarily neighborhoods, and they are not required to form a
fundamental base. In particular, the sets are allowed to be
closed, and their boundaries may have non-empty intersections
with $\mathcal{S}$.

For the chosen class of dynamical systems we formulate
Lyapunov-like conditions that allow to specify
forward invariant sets containing Lyapunov-unstable equilibria on
their boundaries. In the classical method of Lyapunov functions
the role of a Lyapunov function is to assure that an invariant
set, e.g. an equilibrium, has a fundamental base of {  forward}
invariant neighborhoods. In our work we use an extension of
this method in which a substitute of a Lyapunov function is
used to demonstrate existence of a family of {   forward}
invariant sets (not necessarily neighborhoods) associated with
the equilibrium.

The method we use for determining positive invariance of an
individual set is similar in spirit to the second method of
Lyapunov \cite{Lyapunov:1892} and its extensions
\cite{Chetaev:1961}, \cite{Zubov:1964}, \cite{Matrosov:1962},
including  equations with differential inclusions \cite{Barbashin:1961},
\cite{Kolmogorov:1936}, \cite{Nagumo:1942}, \cite{Aubin:1991},
\cite{SICON:2011:Aubin}. Namely, we are looking for  closed
sets containing the origin such that on the boundaries of
these sets the vector-fields in the right-hand side of
(\ref{eq:problem}) are pointing inwards or vanishing. Following
this intuition we demonstrate that there is a set of simple
algebraic conditions, very similar to the ones in the second
method of Lyapunov, enabling us to characterize asymptotic
behavior of solutions for systems with unstable invariant sets.
In particular, these results allow to estimate the domains
of initial conditions, as functions of system parameters, which
are associated with bounded solutions in forward time without
the need to require information about the convergence rate of the
stable part of (\ref{eq:problem}). Parameters of these systems
are not required to be known precisely, and input-output gains
of the systems need not  be defined. Furthermore, in contrast
to our previous results on the same topic
\cite{SIAM_non_uniform_attractivity}, the present conditions allow
to specify domains of initial conditions that lead to
solutions necessarily escaping from a neighborhood of the
equilibria in question.

The paper is organized as follows. In Section
\ref{sec:Formulation} we formulate the problem and specify the main
assumptions. Section \ref{sec:Results} contains the main
results of the paper.

Section \ref{sec:Examples} presents illustrative examples
showing how the results can be applied to 1) derive estimates
of attractor basins for
(\ref{eq:simple_system:1})--(\ref{eq:simple_system:3}), 2)
solve the phase synchronization problem described in Example
\ref{example:0}, and 3) to design an adaptive control scheme
for a class of systems with general nonlinear parametrization.
Section \ref{sec:Conclusion} concludes the paper. Auxiliary
technical results are presented in the Appendix (Section
\ref{sec:Appendix}).

\section{Problem Formulation}\label{sec:Formulation}
Consider system (\ref{eq:problem})
\[
\begin{split}
\dot{x}&=f(x,\lambda,t),\\
\dot{\lambda}&=g(x,\lambda,t),
\end{split} \tag{\ref{eq:problem}}
\]
where the vector-fields $f: \
\Real^n\times\Real\times\Real\rightarrow\Real^n$, $g: \
\Real^n\times\Real\times\Real\rightarrow\Real$ are continuous
and locally Lipschitz w.r.t. $x$, $\lambda$ uniformly in $t$.
 Recall that the point $x=0,
\lambda=0$ is an equilibrium of (\ref{eq:problem}), that $x=0$ is a weak attractor for (\ref{eq:problemFroz}) at $\lambda=0$, and that $\underline{\omega}(p)$, $p>0$ is the set which is forward invariant for all $\lambda\in[0,p]$ w.r.t. the dynamics of (\ref{eq:problemFroz}).

Let  $\mathcal{D}$ be an open subset of $\Real^n$  and
${\Lambda}=[c_1,c_2], \ c_1\leq 0, \ c_2 > 0$, be
an interval. Suppose that the closure
$\overline{\mathcal{D}}$ of $\mathcal{D}$ contains the
origin, and denote
$\mathcal{D}_{\Omega}=\overline{\mathcal{D}}\times{\Lambda}\times\Real$.
Finally,  we suppose that the right-hand side of
(\ref{eq:problem}) satisfies  Assumptions \ref{assume:stable},
\ref{assume:unstable} below.
\begin{assume}\label{assume:stable}There exists a function $V:\Real^n\rightarrow \Real$,  $V\in\mathcal{C}^{0}$, differentiable everywhere except possibly at the origin, and five functions of one variable, $\underline{\alpha},\bar{\alpha}\in\mathcal{K}_{\infty}$, $\alpha:\Real_{\geq 0}\rightarrow\Real$,
$\alpha\in\mathcal{C}^{0}([0,\infty))$, $\alpha(0)=0$, $\beta: \ \Real_{\geq 0}\rightarrow\Real_{\geq 0}$, $\beta\in\mathcal{C}^{0}([0,\infty))$, $\varphi\in\mathcal{K}_0$
such that  for every   $(x,\lambda,t)\in{(\overline{\mathcal{D}}\setminus\{0\})\times\Lambda\times\Real}$  the following properties hold:
\begin{equation}\label{eq:stable}
\begin{split}
\underline{\alpha}(\|x\|)\leq V(x) \leq \bar{\alpha}(\|x\|), \ \ \frac{\pd V}{\pd x}f(x,\lambda,t)\leq  \alpha(V(x)) + \beta(V(x))\varphi(|\lambda|).
 \end{split}
\end{equation}
\end{assume}
Assumption \ref{assume:stable} holds, for example, for systems
in which the term $({\pd V}/{\pd x})f(x,\lambda,t)$ can be
bounded from above as follows: there exist $\alpha_0,\beta_0,\varphi\in\mathcal{K}$ such that
\begin{equation}\label{eq:stable:2:0}
\frac{\pd V}{\pd x}f(x,\lambda,t)\leq - \alpha_0(\|x\|)+\beta_0(\|x\|)\varphi(|\lambda|)  \ \mbox{for \ all} \ (x,\lambda,t)\in{(\overline{\mathcal{D}}\setminus\{0\})\times\Lambda\times\Real}.
\end{equation}
Indeed, (\ref{eq:stable}) follows immediately  from
(\ref{eq:stable:2:0}) with
$\alpha(\cdot)=-\alpha_0\circ\bar{\alpha}^{-1}(\cdot)$,
$\beta(\cdot)=\beta_0\circ\underline{\alpha}^{-1}(\cdot)$, and
$\underline{\alpha}^{-1}(\cdot)$, $\bar{\alpha}^{-1}(\cdot)$
is the inverse of $\underline{\alpha}(\cdot)$,
$\bar{\alpha}(\cdot)$ respectively. In this case Assumption
\ref{assume:stable} states  that the zero solution of (\ref{eq:problemFroz}) at $\lambda=0$ is
 globally asymptotically stable in the sense of Lyapunov,
and  $V(\cdot)$ is the corresponding Lyapunov function.
Notice, however, that Lyapunov stability of the zero solution of
(\ref{eq:problemFroz}) at $\lambda=0$ is not needed for the assumption to hold. System
(\ref{eq:simple_system:2}) is an example of a system in which
the origin is unstable equilibrium and yet Assumption
\ref{assume:stable}  is satisfied with $V(x)=x^2$, $\mathcal{D}=\{x| \ x\in\Real_{>0}\}$  (see Section \ref{sec:examples:academic} for more details).   Finally, we
remark that despite the right-hand side of (\ref{eq:problem})
is allowed to be time-varying, we restrict our consideration to
systems for which the function $V(\cdot)$ does not depend on time
explicitly.

Let us now proceed with detailing  the
requirements for the function $g(\cdot,\cdot,\cdot)$. These are
presented in Assumption \ref{assume:unstable}.
\begin{assume}\label{assume:unstable}  There exist functions $\delta, \xi \in\mathcal{K}_0$ such that the following inequality holds for all $(x,\lambda,t)\in{\mathcal{D}_{\Omega}}$:
\begin{equation}\label{eq:g:Lipschitz:monotone:1}
\begin{split}
- \xi (|\lambda|) - \delta(\|x\|) &\leq  g(x,\lambda,t)\leq 0.
\end{split}
\end{equation}
\end{assume}
Assumption \ref{assume:unstable} reflects the fact that
derivative $\dot{\lambda}$ does not change sign for all
$(x,\lambda,t)\in \mathcal{D}_{\Omega}$. Without loss of
generality we consider the case when $\lambda$ is
non-increasing with time. Alternative formulations of our
conclusions for the case when (\ref{eq:g:Lipschitz:monotone:1})
is replaced with $0\leq  g(x,\lambda,t)\leq  \delta(\|x\|) +
\xi(|\lambda|)$, are readily available. (In
this case one may also need to redefine $\Lambda$
as an interval $[c_1,c_2]$, $c_1<0$, $c_2\geq 0$.)

We aim to formulate a list of conditions that would allow us  to
estimate forward invariant sets of (\ref{eq:problem})
and, specifically, those in which the solutions of
(\ref{eq:problem}) remain bounded. These conditions  are
provided in the next section.

\section{Main Results}\label{sec:Results}

Before  providing formal statements of the results let us
briefly comment on the internal structure of the section. We
begin with Section \ref{subsec:Positive_invariance} presenting
conditions for the existence of forward invariant sets
for (\ref{eq:problem}) containing non-trivial bounded solutions in
forward time. The conditions  are constructive, i.e. not only
existence of such sets is guaranteed but also their
boundaries are explicitly provided. Two alternative statements
of the results are discussed: one is limited to the case of
differentiable boundaries (Lemma
\ref{lem:boundedness_unstable}, Section \ref{subsec:Positive_invariance}), and the other being applicable
to non-differentiable boundaries (Lemma
\ref{lem:boundedness_unstable_2}, Section \ref{subsec:Positive_invariance}). Estimates of the sets corresponding to solutions escaping the origin are provided
in Lemmas \ref{lem:unboundedness_unstable}, \ref{lem:unboundedness_unstable_2} in Section \ref{subsec:Positive_invariance:escaping}.

\subsection{Forward Invariance}\label{subsec:Positive_invariance}

Our first result is provided in the lemma below.
\begin{lem}[Boundedness 1]\label{lem:boundedness_unstable} Let system (\ref{eq:problem}) be given and satisfy Assumptions \ref{assume:stable}, \ref{assume:unstable}.
Suppose that \begin{itemize}
\item[(C1)]  there exist  a function $\psi: \
    \psi\in\mathcal{K}\cap\mathcal{C}^{1}((0,\infty))$
and $a\in\Real_{>0}$ such that for all $V\in(0,a]$
\begin{equation}\label{eq:positive:invariance:attracting}
\frac{\pd \psi(V)}{\pd V}\left[\alpha(V)+\beta(V)\varphi(\psi(V))\right]+\delta\left(\underline{\alpha}^{-1}(V)\right)+\xi\left(\psi(V)\right)\leq 0,
\end{equation}
\end{itemize}
and suppose, in addition to C1
\begin{itemize}
\item[(C2)]   the set $\underline{\omega}(\psi(a))$ exists, and either the set
    $\overline{\mathcal{D}}$  contains
    { $\underline{\omega}(\psi(a))$}, or the ball $\{x \ | \ x\in\Real^n, \ \|x\|\leq
    \underline{\alpha}^{-1}(a)\}$ is in $\mathcal{D}$;
\item[(C3)] the set $\Omega_a\smallsetminus\{(0,0)\}$, where
\begin{equation}\label{eq:positive:invariance:domain}
\Omega_a=\{(x,\lambda) \ | \ x\in{ \underline{\omega}(\psi(a))}, \ \lambda\in\Real_{\geq 0}, \ \psi(a)\geq\lambda\geq\psi(V(x)), \ V(x)\in[0,a]\}
\end{equation}
is contained in the interior of
$\overline{\mathcal{D}}\times \Lambda$.
\end{itemize}
Then
\begin{itemize}
\item[(a) ]$\Omega_a$  is forward invariant with respect to
    (\ref{eq:problem}); that is: solutions of
    (\ref{eq:problem})   starting in $\Omega_a$  at
    $t=t_0$ are defined for all $t\geq t_0$, remain in
    $\Omega_a$, and are  bounded.
\end{itemize}
Furthermore, for every solution of (\ref{eq:problem}) starting in $\Omega_a$
\begin{itemize}
\item[(b) ] there exists a limit
\begin{equation}\label{eq:positive:invariance:limit:1}
  \lim_{t\rightarrow\infty}\lambda(t)=\lambda', \;\; \lambda'\in [0,\psi(a)] .
\end{equation}
\item[(c)] If, in addition, the function $g(x,\lambda,
    \cdot )$ is uniformly continuous then:
\begin{equation}\label{eq:positive:invariance:limit:1.5}
\lim_{t\rightarrow\infty}g(x(t),\lambda',t)=0.
\end{equation}
\end{itemize}
\end{lem}

\begin{figure}
\centering
\includegraphics[width=\textwidth]{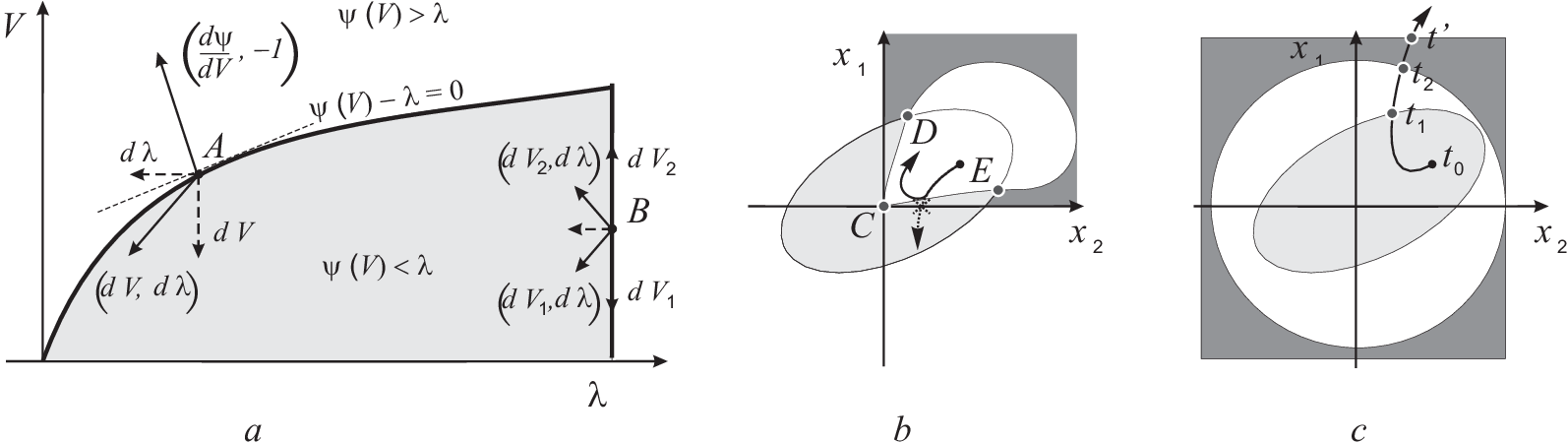}
\caption{Geometric interpretation of conditions of Lemma
\ref{lem:boundedness_unstable}. Panel $a$: condition C1.  Panel $b$: condition C2, first
alternative. Panel $c$: condition C2, second alternative. A more detailed explanation is provided in the text.}\label{fig:positive:invariance:attracting}
\end{figure}

Before providing a proof of the lemma let us first comment
on its conditions. Condition C1 is the actual criterion of forward
invariance. Let us suppose, for simplicity, that
$\overline{\mathcal{D}}$ coincides with $\Real^n$,    and conditions C1--C3 are satisfied for some positive $a$. In this case, similar to the classical  inequality in the
method of Lyapunov functions for systems of ordinary
differential equations $\dot{x}=F(x,t)$, $F\in\mathcal{C}^{0}(\Real^n\times\Real)$:
$(\pd V/\pd x) F(x,t)\leq 0$, inequality
(\ref{eq:positive:invariance:attracting}) guarantees that all
solutions of (\ref{eq:problem}) starting in
(\ref{eq:positive:invariance:domain}) at $t=t_0$ exist for all
$t\geq t_0$ and are bounded in forward time. Condition C2  is a
sort of domestication requirement. It is used in the proof to
ensure that every solution of (\ref{eq:problem}) starting in
$\Omega_a$ and leaving $\Omega_a$ through the boundary
$\lambda=\psi(V(x))$  must necessarily contain a segment
intersecting the boundary $\lambda=\psi(V(x))$ and laying
entirely in $\mathcal{D}\times(0,\psi(a)]$. Condition C3 is a
technical requirement ensuring that every solution crossing
through the boundary $\lambda=\psi(V(x))$, $\lambda\neq 0$ of
$\Omega_a$ (if, of course, such a solution exists) at $t=t'$
will remain in
$\overline{\mathcal{D}}\times\Lambda$ over a
non-empty interval $[t',t'']$, $t''>t'$. A geometric
interpretation of these conditions is provided in Fig.
\ref{fig:positive:invariance:attracting}.

Condition C1 is illustrated in {\it panel} $a$. The vector
$(\pd \psi/\pd V, -1)$ is normal to the curve $\lambda =
\psi(V)$ at the point $A$. Since $\pd \psi/\pd V >0$, it is
always pointing in the direction of  $\lambda<\psi(V)$. One can
easily see that (\ref{eq:positive:invariance:attracting})
implies that for all $V\in(0,a]$ $\left.\frac{\pd \psi}{\pd V} \dot{V} -
\dot{\lambda}\right|_{\lambda=\psi(V)}\leq 0$. Hence, according to
(\ref{eq:positive:invariance:attracting}), the vector
$(\dot{V},\dot{\lambda})$ is pointing in the direction of
$\lambda\geq\psi(V)$ on the surface $\lambda=\psi(V)$.
Condition C2 is illustrated in panels $b$ and $c$. {\it Panel}
$b$ shows the case when the set $\overline{\mathcal{D}}$
contains { $\underline{\omega}(\psi(a))$}. Set $\mathcal{D}$ is the
dark grey rectangle, the set $\lambda\geq \psi(V(x))$ for some
fixed value of $\lambda$ is depicted as a light grey ellipse,
and { $\underline{\omega}(\psi(a))$} is the white area. The condition
states that only those crossings through the boundary
$\lambda=\psi(V(x))$ (segment DE in the figure) are allowed
which occur in the white area. Solutions cannot cross segments
CD and CE and hence must remain in $\overline{\mathcal{D}}$.
{\it Panel} $c$ concerns the second alternative, i.e. when the
ball $\{x \ | x\in\Real^n, \ \|x\|\leq \underline{\alpha}^{-1}(a)\}$ is in
$\mathcal{D}$.  The white area depicts the ball  $\{ x \ | \ x\in\Real^n, \ \|x\|\leq
\underline{\alpha}^{-1}(a)\}$. The ball contains sets
$\{x\ | \ x\in\Real^n, \   \ V(x)\leq c, \ c\in[0,a]\}$ as subsets.
Given that the function $\psi(\cdot)$ is non-decreasing and strictly
monotone, it is clear that the ball contains   $\{ x \ |
x\in\Real^n, \   \  \lambda\geq \psi(V(x)),  \
\lambda\in[0,\psi(a)]\}$.  The condition therefore reflects
that any solution  $x(\cdot;t_0,x_0,\lambda_0)$ leaving the set
$\overline{D}$ at $t=t'$ must necessarily cross through the
surfaces $\lambda=\psi(V(x))$, $\lambda\in[0,\psi(a)]$ and
$\|x\|=\underline{\alpha}^{-1}(a)$ at $t=t_1$ and  $t=t_2$,
$t_2\geq t_1$ respectively.

{\it Proof of Lemma \ref{lem:boundedness_unstable}.}
{The proof of the lemma is split into two parts.
In the first part we show that conditions C2, C3, and the fact that Assumptions
\ref{assume:stable}, \ref{assume:unstable} hold guarantee that every solution of (\ref{eq:problem}) satisfying the initial condition $x(t_0)=x_0$, $\lambda(t_0)=\lambda_0$, $(x_0,\lambda_0)\in\Omega_a$
 must either 1)
remain in ${ \underline{\omega}(\psi(a))} \times[0,\psi(a)]$ for $t\geq t_0$ (and
consequently in $\overline{\mathcal{D}}\times[0,\psi(a)]$) as long as $\lambda(t;t_0,x_0,\lambda_0)\geq 0$
(first alternative), or 2) if it leaves the set
$\overline{\mathcal{D}}\times[0,\psi(a)]$ at some $t\geq t_0$
then it should first cross the boundary $\lambda=\psi(V(x))$, $\lambda\in(0,\psi(a)]$
(second alternative) in
$\overline{\mathcal{D}}\times[0,\psi(a)]$. This ensures that
inequalities (\ref{eq:stable}),
(\ref{eq:g:Lipschitz:monotone:1}) in Assumptions
\ref{assume:stable}, \ref{assume:unstable} must hold along the
solutions of  (\ref{eq:problem}) starting in $\Omega_a$ at
$t_0$ for $t\geq t_0$ as long as they remain in
$\Omega_a\smallsetminus\{(0,0)\}$.  Furthermore, if the solution
crosses through the boundary of the set
$\Omega_a\smallsetminus\{(0,0)\}$ at some $t\geq t_0$ then it must
necessarily satisfy Assumptions \ref{assume:stable},
\ref{assume:unstable} over a non-empty interval $[t,t']$,
$t'>t$ because the set $\Omega_a\smallsetminus\{(0,0)\}$ is in the
interior of
$\overline{\mathcal{D}}\times \Lambda$.

In the second part
of the proof we use this property to show that condition
(\ref{eq:positive:invariance:attracting}) is incompatible with
the assertion that solutions  of (\ref{eq:problem}) startingh $\Omega_a$ at $t=t_0$ may intersect the boundary
$\lambda=\psi(V(x))$  at $t\geq t_0$.

{\it  Part 1.}   Let   $(x_0,\lambda_0)$ be a point in $\Omega_a$. It is clear that solutions of (\ref{eq:problem}) exist at
least locally and   are unique.
 According to the first alternative of condition C2, that   $\overline{D}\supset { \underline{\omega}(\psi(a))}$, components $x(t;t_0,x_0,\lambda_0)$,   $x_0\in { \underline{\omega}(\psi(a))}$, $\lambda_0\in[0,\psi(a)]$ of the solutions of (\ref{eq:problem}) must belong to $\overline{\mathcal{D}}$ for $t\geq t_0$ as long as
$\lambda(t;t_0,x_0,\lambda_0)\geq 0$ (see Remark \ref{rmk:ExtInv}).


 Consider the second alternative of C2;   $\mathcal{D}$ contains the ball $\{ x \ | \ x\in\Real^n, \  \|x\|\leq
\underline{\alpha}^{-1}(a)\}$.   Since the right-hand side of (\ref{eq:problem})
is locally Lipschitz, the equilibrium solution $x(\cdot;t_0,0,0)\equiv 0$, $\lambda(\cdot;t_0,0,0)\equiv0$ is unique. Thus solutions $x(\cdot;t_0,x_0,\lambda_0)$, $\lambda(\cdot;t_0,x_0,\lambda_0)$,  $(x_0,\lambda_0)\in\Omega_a$ cannot escape the domain $\mathcal{D}\times[0,\psi(a)]$ through the point $(0,0)$. Let us show that
if there is a solution $x(\cdot;t_0,x_0,\lambda_0)$, $\lambda(\cdot;t_0,x_0,\lambda_0)$, $(x_0,\lambda_0)\in\Omega_a$ of (\ref{eq:problem}) that is leaving the set $\overline{\mathcal{D}}\times[0,\psi(a)]$ at some $t>t_0$ then it must first cross
the boundary $\lambda=\psi(V(x))$, $\lambda\in(0,\psi(a)]$.
Let this
not be the case and suppose that there exists a solution of
(\ref{eq:problem}) touching the boundary of
$\overline{\mathcal{D}}$ without crossing through
$\lambda=\psi(V(x))$, $\lambda\in(0,\psi(a)]$. This means that
there exists $t'> t_0$  such that
\begin{subequations}
\begin{eqnarray}\label{eq:lem:1:second_contradiction}
\lambda(t';t_0,x_0,\lambda_0)&\geq& \psi(V(x(t';t_0,x_0,\lambda_0)), \\ \lambda(t';t_0,x_0,\lambda_0)&\in&(0, \psi(a)],\label{eq:lem:1:second_contradiction:2}
\end{eqnarray}
\end{subequations}
i.e. no crossing occurred, and yet the   point
 $x(t';t_0,x_0,\lambda_0)$ is on the boundary of
 $\overline{\mathcal{D}}$. It is therefore clear that the
following must hold:
$\psi(\underline{\alpha}(\|x(t';t_0,x_0,\lambda_0)\|)) > \psi(a)$.
On the other hand, according to Assumption \ref{assume:stable},
we have that $\psi(\underline{\alpha}(\|x(t';t_0,x_0,\lambda_0)\|))\leq
\psi(V(x(t';t_0,x_0,\lambda_0)))$, and hence  $\psi(V(x(t';t_0,x_0,\lambda_0)))>\psi(a)$.
The latter inequality together with
(\ref{eq:lem:1:second_contradiction}) result in
$\lambda(t';t_0,x_0,\lambda_0)>\psi(a)$. This, however, contradicts
to (\ref{eq:lem:1:second_contradiction:2}).

{\it  Part 2.} We claim that any solution of (\ref{eq:problem})
passing through $(x_0,\lambda_0)$ at $t_0$ is defined for all
$t\geq t_0$ and remains in $\Omega_a$ for all $t\geq t_0$. Let
us first demonstrate that solutions of (\ref{eq:problem})
starting in $\Omega_a\smallsetminus\{(0,0)\}$  cannot leave the set
through the boundary $\lambda=\psi(V(x))$, $V(x)\in(0,a]$. Assume that this is not the case. Pick an arbitrary point
$(x_0,\lambda_0)\in\Omega_a\smallsetminus\{(0,0)\}$, and let
$\phi(\cdot;t_0,x_0,\lambda_0)=(x(\cdot;t_0,x_0,\lambda_0),\lambda(\cdot;t_0,x_0,\lambda_0))$, $(x_0,\lambda_0)\in\Omega_a\smallsetminus\{(0,0)\}$
be the   maximal solution of (\ref{eq:problem}), and let $\mathcal{T}=[t_0,t_{\max})$ be the interval of its definition for $t>t_0$.      Suppose that $\phi(\cdot;t_0,x_0,\lambda_0)$
can cross through the boundary, i.e. there exist
$t'\in\mathcal{T}$ such that
$\lambda(t';t_0,x_0,\lambda_0)<\psi(V(x(t';t_0,x_0,\lambda_0)))$.
Condition C3 states that $\Omega_a$ is in the
interior of
$\overline{\mathcal{D}}\times\Lambda$. Hence
without loss of generality we can suppose that
$\phi(t';t_0,x_0,\lambda_0)\in
\overline{\mathcal{D}}\times\Lambda$.

Consider the function $p:\mathcal{T}\rightarrow\Real$,
$p(t)=\psi(V(x(t;t_0,x_0,\lambda_0)))-\lambda(t;t_0,x_0,\lambda_0)$.
The function   $p(\cdot)$   is continuous in $\mathcal{T}$. Thus there is
a non-empty interval $[t_1,t']\subset\mathcal{T}$, such that
$p(t_1)=0$ and $p(t)>0$ for all  $t\in(t_1,t']$. Moreover, given that
$\psi\in\mathcal{C}^{1}((0,\infty))$, for   every $t\in(t_1,t']$ the derivative $\dot{p}(t)$
exists, and   $\dot{p}(\cdot)$   is locally bounded in $\mathcal{T}$.   Therefore $p(\cdot)$   is absolutely continuous in $[t_1,t']$, and
$p(t)=\int_{t_1}^{t} \dot{p}(\tau)d\tau > 0, \ t\in(t_1,t']$.
According to the mean-value theorem there exists a
$\tau\in(t_1,t]$ such that
\[
\begin{split}
p(t)&= (t-t_1) \dot{p}(\tau)= (t-t_1)\big[\frac{\pd \psi}{\pd  V}\frac{\pd V}{\pd x}f(x(\tau;t_0,x_0,\lambda_0),\lambda(\tau;t_0,x_0,\lambda_0),\tau)\\
&-g(x(\tau;t_0,x_0,\lambda_0),\lambda(\tau;t_0,x_0,\lambda_0),\tau)\big]
\end{split}
\]
Using the  fact that the function $\psi(\cdot)$ is non-decreasing,
i.e. $\pd \psi/\pd V\geq 0$ for all $V\in(0,a]$, and invoking
Assumptions \ref{assume:stable}, \ref{assume:unstable} we
derive that
\[
\begin{split}
 \dot{p}(\tau)&\leq \frac{\pd \psi}{\pd  V} \left[\alpha(V(x(\tau;t_0,x_0,\lambda_0))) + \beta(V(x(\tau;t_0,x_0,\lambda_0)))\varphi(\lambda(\tau;t_0,x_0,\lambda_0))\right]\\
& +  \xi (\lambda(\tau;t_0,x_0,\lambda_0)) + \delta(\|x(\tau;t_0,x_0,\lambda_0)\|).
\end{split}
\]
The functions $\xi(\cdot),\delta(\cdot),\varphi(\cdot)$ are
non-decreasing, and
$\lambda(\tau;t_0,x_0,\lambda_0)$ $<$ $\psi(V(x(\tau;t_0,x_0,\lambda_0)))$. Hence
 invoking condition C1 of the lemma, we can conclude that:
\begin{equation}\label{eq:lem:1:contradiction}
\begin{split}
 \dot{p}(\tau)&\leq  \frac{\pd \psi}{\pd  V} \left[\alpha(V(x(\tau;t_0,x_0,\lambda_0))) + \beta(V(x(\tau;t_0,x_0,\lambda_0)))\varphi(\psi(V(x(\tau;t_0,x_0,\lambda_0))))\right]\\
& +  \xi (\psi(V(x(\tau;t_0,x_0,\lambda_0)))) + \delta(\underline{\alpha}^{-1}(V(x(\tau;t_0,x_0,\lambda_0))))\leq 0.
\end{split}
\end{equation}
On the other hand, since $p(t)>0$  and $(t-t_1)>0$ for all $t\in(t',t_1]$,  the
following must hold: $\dot{p}(\tau)>0$. This contradicts
(\ref{eq:lem:1:contradiction}), and hence the statement that
the solution $\phi(\cdot;t_0,x_0,\lambda_0)$
crosses the boundary $\lambda=\psi(V(x))$ at some $t'\in\mathcal{T}$ in
$\overline{\mathcal{D}}\times \Lambda\smallsetminus\{(0,0)\}$
in finite time is not true.

It is also clear that $\phi(\cdot;t_0,x_0,\lambda_0)$ cannot escape
$\Omega_a$ through the boundary $\lambda=\psi(a)$ at any $t\in\mathcal{T}$ since the
derivative $\dot{\lambda}(t)$ is non-positive for all $t\in\mathcal{T}$. Finally, notice
that the only remaining subset of the boundary of $\Omega_a$
through which the solutions may escape is the set $\{
(x,\lambda) \ | \ x\in\Real^n, \ \lambda\in\Real, \ V(x)=0,\
\lambda=0\}$. This set, however, is the equilibrium of
(\ref{eq:problem}),  and the equilibrium solution $\phi(\cdot;t_0,0,0)=0$ is unique.

Thus $\phi(t;t_0,x_0,\lambda_0)\in\Omega_a$ for all $t\in\mathcal{T}$.  Noticing that the set
$\Omega_a$ is bounded, and that the right-hand side of
(\ref{eq:problem}) is locally Lipschitz we conclude that  $\phi(\cdot;t_0,x_0,\lambda_0)$ is defined on $[t_0,\infty)$
and bounded.   Given that $(x_0,\lambda_0)$ was an arbitrary
point of $\Omega_a\smallsetminus\{(0,0)\}$, that  $t_0$ was chosen
arbitrary,  and that  the origin is the equilibrium of
(\ref{eq:problem}), we conclude that all solutions of
(\ref{eq:problem}) passing through $\Omega_a$ remain in
$\Omega_a$ in forward time.

Property (\ref{eq:positive:invariance:limit:1}) is an immediate
consequence of the Bolzano-Weierstrass theorem. In order to see
that property (\ref{eq:positive:invariance:limit:1.5}) holds, we
notice that the integral
$\lambda(t;t_0,x_0,\lambda_0)=\lambda_0+\int_{t_0}^t
g(x(\tau;t_0,x_0,\lambda_0),\lambda(\tau;t_0,x_0,\lambda_0),\tau)d\tau$
converges. Taking into account boundedness of  $x(\cdot;t_0,x_0,\lambda_0)$, $\lambda(\cdot;t_0,x_0,\lambda_0)$ on $[t_0,\infty)$,
the fact that   $g(\cdot,\cdot,t)$ is locally Lipschitz  uniformly in $t$   and is uniformly continuous in $t$ we conclude
that
$g(x(\cdot;t_0,x_0,\lambda_0),\lambda(\cdot;t_0,x_0,\lambda_0),\cdot)$, is uniformly continuous on $[t_0,\infty)$.   Thus, according
to the Barbalat's lemma (see e.g. \cite{Bermant:1963}),
$$\lim_{t\rightarrow\infty}g(x(t;t_0,x_0,\lambda_0),\lambda(t;t_0,x_0,\lambda_0),t)=0.$$
Hence the result follows as a consequence of
(\ref{eq:positive:invariance:limit:1}) and continuity of
$g(x,\cdot,t)$ in $\Lambda$. $\square$ }

\begin{rmk}\label{rem:measure} \normalfont Notice that if the interior of ${ \underline{\omega}(\psi(a))}$
is not empty and its closure contains the origin then the measure of $\Omega_a$, $a>0$ is not zero.  Indeed, since the function $\psi(\cdot)$ is strictly monotone, $\psi(a)>0$. Pick a $\lambda'\in (0,\psi(a))$, and let $\mathcal{B}(x,r)$ denote an open ball of radius $r$ in $\Real^n$  centered at $x$. The ball $\mathcal{B}(0,(\overline{\alpha}^{-1}\circ \psi^{-1})(\lambda'))$ is contained in
the set  $\{x\in\Real^n | \ x: \ \psi^{-1}(\lambda')> V(x)\}$. Since the closure of
${ \underline{\omega}(\psi(a))}$ contains the origin and interior of ${ \underline{\omega}(\psi(a))}$
is open we can conclude that there exists a ball
$\mathcal{B}(x',r)$  such that $\mathcal{B}(x',r)\subset
{ \underline{\omega}(\psi(a))}$ and $\mathcal{B}(x',r)\subset
\mathcal{B}(0,(\overline{\alpha}^{-1}\circ
\psi^{-1})(\lambda'))$. Thus the set
$\mathcal{B}(x',r)\times\{\lambda'\}$ is in the interior of
$\Omega_a$. It is also clear that the set
$\overline{\mathcal{B}}(x',r)\times[\lambda',\lambda'']$,
$\psi(a)>\lambda''>\lambda'$
 is in the interior of $\Omega_a$ provided that $\lambda''$ is sufficiently close to $\lambda'$. Thus the measure of $\Omega_a$ is not zero.
\end{rmk}

\begin{rmk}\label{rem:simple_necessary_condition} \normalfont According to the assumptions of Lemma \ref{lem:boundedness_unstable},  $d\psi(V(x(t;x_0,\lambda_0)))/dt$ is non-positive in the set $\{(\lambda,x)  \ \lambda\in\Real_{>0}, \ x\in\Real^n |  \ \lambda=\psi(V(x)), \  V(x)\in[0,a]\}$.
Hence if $V(x)\in[0,a]$ and $t\geq t_0$ then
\begin{equation}\label{eq:necessary:lemma_1}
\frac{\pd V}{\pd x}f(x,\psi(V(x)),t)\leq  0.
\end{equation}
Therefore, subject to the assumptions of the lemma,
(\ref{eq:necessary:lemma_1}) may be used as a necessary
condition for testing positive invariance of the sets
$\Omega_a$.
\end{rmk}
\begin{rmk} \label{rem:simple_generalization} \normalfont If the second alternative of C2 holds then the requirement that $\underline{\omega}(\psi(a))$ exists is not necessary, and the definition of $\Omega_a$ may be replaced with: $\Omega_a=\{(x,\lambda) \ | \ x\in \Real^n, \ \lambda\in\Real_{\geq 0}, \ \psi(a)\geq\lambda\geq\psi(V(x)), \ V(x)\in[0,a]\}$.
\end{rmk}
 Notice also, that if  $\mathcal{D}$,
$\Lambda$ coincide with $\Real^n$ and $\Real$,
respectively, then Lemma \ref{lem:boundedness_unstable} reduces
to the much simpler statement below.

\begin{cor}\label{cor:boundedness_unstable} Consider  system (\ref{eq:problem}), and let $\mathcal{D}=\Real^n$, $\Lambda=\Real$. Suppose that Assumptions \ref{assume:stable}, \ref{assume:unstable} hold, there exists  a function
$\psi: \ \psi\in\mathcal{K}\cap\mathcal{C}^{1}((0,\infty))$ and a
positive constant $a$ such that
(\ref{eq:positive:invariance:attracting}) holds. Then the set
\begin{equation}\label{eq:positive:invariance:domain:cor}
\Omega_a=\{(x,\lambda) \ | \  x\in\Real^n, \ \lambda\in\Real_{\geq 0}, \ \psi(a)\geq\lambda\geq\psi(V(x)), \ V(x)\in[0,a]\}
\end{equation}
is forward  invariant and conclusions
(b), (c) of Lemma \ref{lem:boundedness_unstable} hold.
\end{cor}

Lemma \ref{lem:boundedness_unstable} requires that the function
$\psi(\cdot)$ used in the definition of the set $\Omega_a$ is
strictly increasing and differentiable. Occasionally, a
need might arise for functions $\psi(\cdot)$ which are not
differentiable or not strictly monotone. While the requirement
of strict monotonicity can be traded for the weaker constraint
that the function $\psi(\cdot)$ is non-decreasing, without
significant alternations to the statements of Lemma
\ref{lem:boundedness_unstable} and Corollary
\ref{cor:boundedness_unstable} and their proofs, dealing with
the issue of differentiability involves replacing
(\ref{eq:positive:invariance:attracting}) with a set of
different invariance conditions. In what follows we present
these modified conditions involving the notion of a star shaped
set. A formal definition and basic properties of star shaped
sets and functions are provided in Appendix. A brief summary
of notions that  are essential for the formulation of the
results is provided below.

Recall that a set $\Omega\subset\Real^n$ is called
{\it star-shaped with respect to the origin} if every segment
connecting the origin with a  point $p\in\Omega$ lies entirely
in $\Omega$. Clearly, convexity of a set containing the origin
implies that it is also star-shaped with respect the origin
(since a convex set is star-shaped w.r.t. every point of that
set). The star shaped envelop of a set $D$ (w.r.t. the origin)
is the minimal star shaped set (w.r.t. the origin) including
$D$; that is, every star shaped set (w.r.t. the origin)
including $D$ as a subset must necessarily include the star shaped
envelop (w.r.t. the origin) of $D$ as well. This set is denoted
as $\mathrm{star}(D)$. (In the Appendix, for a star
shaped envelop of a set with respect to a point $x$, the
notation $\mathrm{star}_x(D)$ is used. For notational convenience we
omit the subscript $0$ when refering to star shaped sets,
envelopes, and functions with respect to the origin.)
The {\it epigraph}, respectively {\it hypograph}, of a function $f:\Real^n\rightarrow\Real$,
or simply $\mathrm{epi}(f)$, is the set in $\Real^{n+1}$:
$\mathrm{epi}(f)=\{(x,\mu)| \ x\in\Real^n,\ \mu\in\Real,
f(x)\leq \mu \}$, respectively $\mathrm{hyp}(f)=\{(x,\mu)| \ x\in\Real^n,\ \mu\in\Real,
f(x)\geq \mu \}$. A function $f:\Real^n\rightarrow\Real$ is
called {\it star shaped} with respect to the origin iff its
epigraph is a star shaped set with respect to the origin.
The convex envelope of a function $f:\Real^n\rightarrow \Real$ is denoted as
$\mathrm{conv}(f)(x)$ and the star shaped envelop of $f(\cdot)$ (w.r.t.
the origin) is denoted as $\mathrm{star}(f)(x)$. Let
$f:\Real\rightarrow\Real$, and $[0,a]$, $a>0$ be an interval.
We define
\[
\begin{split}
\mathrm{epi}(f_{[0,a]})&=\{ (x,\mu) \ |\ x\in[0,a],\mu\in\Real, \ f(x)\leq \mu \},\\
\mathrm{hyp}(f_{[0,a]})&=\{ (x,\mu) \ |\ x\in[0,a],\mu\in\Real, \ f(x)\geq \mu \},
\end{split}
\]
the epigraph and hypograph of the restriction of $f$ on the
interval $[0,a]$. Now we are ready to formulate the following:

\begin{lem}[Boundedness 2]\label{lem:boundedness_unstable_2} Let system (\ref{eq:problem}) be given, and let it satisfy Assumptions \ref{assume:stable}, \ref{assume:unstable}. Suppose that
\begin{itemize}
\item[(C4)]  there  exists a function $\psi\in\mathcal{K}$,
    such that for some $a\in\Real_{>0}$ the set
    $\mathrm{epi}(\psi_{[0,a]})$ is star-shaped with
    respect to the origin, and for all $V\in[0,a]$
\begin{equation}\label{eq:positive:invariance:attracting_2}
\psi(V)\left[\alpha(V)+\beta(V)\varphi(\psi(V))\right]+ V \left[\delta\left(\underline{\alpha}^{-1}(V)\right)+\xi\left(\psi(V)\right)\right]\leq 0.
\end{equation}

\end{itemize}
Furthermore, let conditions C2,C3 of Lemma
\ref{lem:boundedness_unstable} hold.

Then $\Omega_a$ is forward invariant with respect to
(\ref{eq:problem}), and moreover conclusions (b), (c) of  Lemma
\ref{lem:boundedness_unstable} apply.
\end{lem}
{\it Proof of Lemma \ref{lem:boundedness_unstable_2}.} {The
proof is similar to that of Lemma
\ref{lem:boundedness_unstable}.

As has already been shown, conditions C2, C3, and Assumptions
\ref{assume:stable}, \ref{assume:unstable} imply that if the solution $x(\cdot;t_0,x_0,\lambda_0)$, $\lambda(\cdot;t_0,x_0,\lambda_0)$,  $(x_0,\lambda_0)\in\Omega_a$
leaves the set
$\Omega_a$ at some $t\geq t_0$ then it must
necessarily satisfy Assumptions \ref{assume:stable},
\ref{assume:unstable} over a non-empty interval $[t,t']$,
$t'>t$. We will show now that condition (\ref{eq:positive:invariance:attracting_2}) is
incompatible with the claim that solutions of the system
starting in the set $\Omega_a$ at $t=t_0$ may leave the set at
some $t'\geq t_0$.

Let $(x_0,\lambda_0)\in\Omega_a$, and
$\phi(\cdot;t_0,x_0,\lambda_0)$ be a solution of (\ref{eq:problem}). The solution exists at
least locally; let $\mathcal{T}=[t_0,t_{\max})$ be its maximal
interval of definition for $t\geq t_0$. First, we show that  if  $(x_0,\lambda_0)\neq (0,0)$ then the solution
$\phi(\cdot;t_0,x_0,\lambda_0)$ does not leave $\Omega_a$ through the
boundary $\lambda=\psi(V(x))$, $V(x)\in(0,a]$ at any
$t\in\mathcal{T}$. Suppose that this is not the case, and there
exists a $t'\in\mathcal{T}$ such that
$\lambda(t';t_0,x_0,\lambda_0)<\psi(V(x(t';t_0,x_0,\lambda_0)))$.
Without loss of generality we can assume that the value of $t'$
is such that
$\lambda(t;t_0,x_0,\lambda_0)\geq\psi(V(x(t;t_0,x_0,\lambda_0)))$ for
all  $t\in[t_0,t_1]$, $t_0\leq t_1$, but
$\lambda(t;t_0,x_0,\lambda_0)<\psi(V(x(t;t_0,x_0,\lambda_0)))$ for
$t\in(t_1,t']$. Further, let $t'$ be close enough to $t_0$, so
that $\lambda(t;t_0,x_0,\lambda_0)\neq 0$ for all $t\in[t_0,t']$.
It is clear that making such choice of $t'$ is always possible
because the value of $\lambda(t_1;t_0,x_0,\lambda_0)$ must necessarily be
positive, and that $\lambda(\cdot;t_0,x_0,\lambda_0)$ is continuous on
$\mathcal{T}$.

Introduce the   function $p:[t_0,t']\rightarrow\Real$: $p(t)={V(x(t;t_0,x_0,\lambda_0))}/{\lambda(t;t_0,x_0,\lambda_0)}$. 
The   function $p(\cdot)$   is defined on $[t_0,t']$ and,
moreover, it is continuous and is continuously differentiable on $[t_0,t']$.
Thus $p(\cdot)$ is absolutely continuous on $[t_0,t']$, and
$p(t)=p(t_1)+\int_{t_1}^{t}\dot{p}(\tau)d\tau$.
Notice that   $p(t)>p(t_1)$ for $t\in(t_1,t']$.   Indeed, if $p(t)\leq p(t_1)$ then the point  $(V(x(t;t_0,x_0,\lambda_0)),\lambda(t;t_0,x_0,\lambda_0))$ would
belong to the set  $\mathrm{epi}(\psi_{[0,a]})$ which
contradicts the   earlier established property that
$\psi(V(x(t;t_0,x_0,\lambda_0)))<\lambda(t;t_0,x_0,\lambda_0))$ for all $t\in(t_1,t']$.

In order to see this consider the point
$(V(x(t_1;t_0,x_0,\lambda_0)),\lambda(t_1;t_0,x_0,\lambda_0))$.
According to the choice of $t_1$, this point is from
$\mathrm{epi}(\psi_{[0,a]})$. Given that
$\mathrm{epi}(\psi_{[0,a]})$ is star shaped (condition C4 of
the lemma), every pair $(v,\lambda)$: $\lambda \geq v/p(t_1)$,
$v\in[0,V(x(t_1;t_0,x_0,\lambda_0)]$ belongs to
$\mathrm{epi}(\psi_{[0,a]})$. Noticing that
$\lambda(t;t_0,x_0,\lambda_0)\leq \lambda(t_1;t_0,x_0,\lambda_0)$ for all
$t\in(t_1,t']$,   we can conclude that   the condition: $p(t)\leq p(t_1)$ for all
$t\geq t_1$   implies that
$V(x(t;t_0,x_0,\lambda_0))$ $\leq$ $V(x(t_1;t_0,x_0,\lambda_0))$ for all
$t\in(t_1,t']$. Thus points
$(V(x(t;t_0,x_0,\lambda_0)),\lambda)$ such that $\lambda$ $\geq$ $
V(x(t;t_0,x_0,\lambda_0))/p(t_1)$ are in
$\mathrm{epi}(\psi_{[0,a]})$. This includes
$(V(x(t;t_0,x_0,\lambda_0),\lambda(t;t_0,x_0,\lambda_0))$ for $t\in(t_1,t']$    provided that
$p(t)=V(x(t;t_0,x_0,\lambda_0))$ $/$ $\lambda(t;t_0,x_0,\lambda_0)\leq
p(t_1)$ for   $t\in(t_1,t']$.

Having derived that $p(t)>p(t_1)$   for all $t\in(t_1,t']$ we therefore arrive at:
$\int_{t_1}^{t'}\dot{p}(s)ds>0$. According to the mean-value theorem, there is a
$\tau\in(t_1,t']$  such that
$\dot{p}(\tau)>0$.
Assumptions \ref{assume:stable},\ref{assume:unstable}, however,
imply that
\[
\begin{split}
&\dot{p}(\tau)= \frac{\dot{V}(x(\tau;t_0,x_0,\lambda_0))\lambda(\tau;t_0,x_0,\lambda_0)-\dot{\lambda}(\tau;t_0,x_0,\lambda_0)V(x(\tau;t_0,x_0,\lambda_0)) }{\lambda^2(\tau;t_0,x_0,\lambda_0)}\leq 0
\end{split}
\]
because condition C4 of the lemma demands that
\[
\begin{split}
& \dot{V}(x(\tau;t_0,x_0,\lambda_0))\lambda(\tau;t_0,x_0,\lambda_0)-\dot{\lambda}(\tau;t_0,x_0,\lambda_0)V(x(\tau;t_0,x_0,\lambda_0))\leq\psi(V(x(\tau;t_0,x_0,\lambda_0)))\\
& \times \big[\alpha(V(x(\tau;t_0,x_0,\lambda_0)))+\beta(V(x(\tau;t_0,x_0,\lambda_0)))\varphi(\psi(V(x(\tau;t_0,x_0,\lambda_0))))\big]+\\
& \ \ \ \ \ \ \ \ \ \ \ \ \ \ \ \ \ \ V \big[\delta\big(\underline{\alpha}^{-1}(V(x(\tau;t_0,x_0,\lambda_0)))\big)+\xi\big(\psi(V(x(\tau;t_0,x_0,\lambda_0)))\big)\big]\leq 0.
\end{split}
\]
This is clearly a contradiction, and hence   the solution
$\phi(\cdot;t_0,x_0,\lambda_0)$ does not cross the boundary
$\lambda=\psi(V(x))$, $V\in(0,a]$ of $\Omega_a$ at any $t\in\mathcal{T}$. The rest of
the proof is identical to that of Lemma
\ref{lem:boundedness_unstable}. } $\square$

\begin{figure}
\begin{center}
\includegraphics[width=250pt]{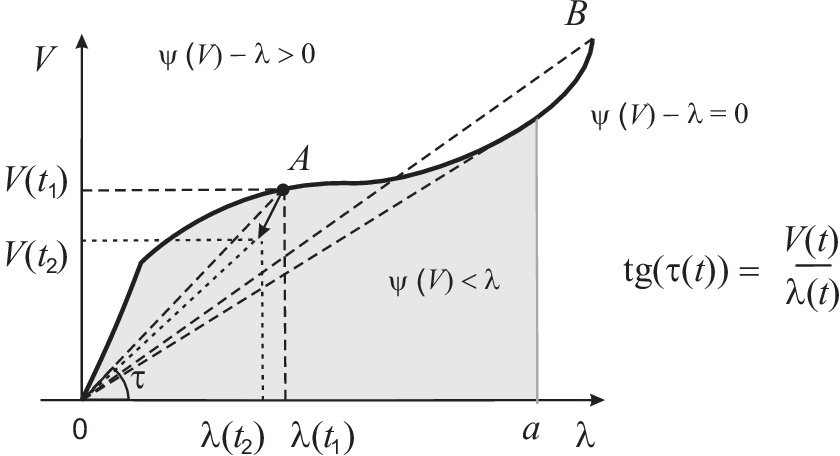}
\end{center}
\caption{Conditions of forward invariance with non-differentiable function $\psi(\cdot)$.
 Star-shapedness of the  $\mathrm{epi}(\psi)$  is needed to guarantee that the vector $(V(t_1)-V(t_2), \lambda(t_1)-\lambda(t_2))$ points in the direction of $\psi(V)<\lambda$ for all points on the curve $\lambda=\psi(V)$, $V\in[0,a]$. Notice that this condition does not hold for the point $B$.}\label{fig:positive:invariance:3}
\end{figure}

 There is a simple geometric interpretation of the conditions of Lemma \ref{lem:boundedness_unstable_2} (see  Fig. \ref{fig:positive:invariance:3}). Consider an interval $(0,a]$ such that the epigraph of $\psi$ for $V\in[0,a]$:
$\mathrm{epi}(\psi_{[0,a]})=\{(V,\lambda) \ | \ V\in [0,a], \
\lambda\in\Real_{\geq 0}, \  \psi(V)< \lambda\}$ is star-shaped
with respect to the origin. It is clear that if the
vector $(\dot{V},\dot{\lambda})$ at the boundary
$\lambda=\psi(V)$, $V\in(0,a]$, is always pointing inside
$\mathrm{epi}(\psi_{[0,a]})$ then
$V(x(t;t_0,x_0,\lambda_0))$, $\lambda(t;x_0,\lambda_0)$,
$(x_0,\lambda_0\in\Omega_a)$ will remain in $\Omega_a$ for all
$t\geq t_0$. A sufficient condition for the latter, as it is
illustrated in  Fig. \ref{fig:positive:invariance:3}, is that the ratio
$p(t)={V(x(t;t_0,x_0,\lambda_0))}/{\lambda(t;t_0,x_0,\lambda_0)}$ is non-increasing with $t$, $t\geq t_0$
for all $V\in[0,a]$, $\lambda=\psi(V)$. Rewriting this
requirement as $\dot{p}(t)\leq 0$,   for all $t\geq t_0$ and invoking the estimates in
Assumptions \ref{assume:stable}, \ref{assume:unstable}  results
in (\ref{eq:positive:invariance:attracting_2}).

\begin{rmk}\label{rmk:V_time_dependent}\normalfont Note that Lemmas \ref{lem:boundedness_unstable}, \ref{lem:boundedness_unstable_2} and Corollary \ref{cor:boundedness_unstable} can be straightforwardly generalized to account for a wider range of systems (\ref{eq:problem}). In particular, in Assumption \ref{assume:stable}, instead of requiring the existence  of  $V:\Real^n\rightarrow\Real$ such that (\ref{eq:stable}) holds one can require that a weaker condition is satisfied: There exists a function of $n+1$ real variables, $V:\Real^n\times\Real\rightarrow\Real$ and two functions of one variables, $\underline{\alpha},\ \bar{\alpha}\in\mathcal{K}_{\infty}$, such that $V(\cdot,\cdot)$ is continuously differentiable except (possibly) on the ray $\{(x,t) \ | \ x\in\Real^n, \ t\in\Real, \ t\geq t_0, \ x=0\}$, and
\begin{equation}\label{eq:stable:general}
\begin{split}
 \underline{\alpha}(\|x\|)\leq V(x,t)\leq \bar{\alpha}(\|x\|), \ \frac{\pd V}{\pd x}f(x,\lambda,t) + \frac{\pd V}{\pd t} \leq  \alpha(V(x,t)) + \beta(V(x,t))\varphi(|\lambda|).
 \end{split}
\end{equation}
It is clear that statements of the results presented so far will remain the same, and also that their proofs will be almost identical under such a modification. The only visible differences would be notational, making technical derivations look more complicated.
\end{rmk}

\subsection{Solutions escaping a neighborhood of the origin}\label{subsec:Positive_invariance:escaping}

In the previous sections we focused predominantly on
developing  boundedness and convergence tests for the solutions
of (\ref{eq:problem}). It may sometimes be desirable to
consider a complementary problem: determine whether a
solution passing through some given point escapes the vicinity
of the origin. As we show below, the very same set of
techniques that was employed in the proofs Lemmas
\ref{lem:boundedness_unstable},
\ref{lem:boundedness_unstable_2}, can be applied for solving
the complementary problem as well.

Naive intuition suggests that criteria specifying domains of
initial conditions corresponding to the solutions escaping
small neighborhoods of the origin could be in a dual relation
to the convergence and boundedness tests we formulated earlier.
In some sense, this is indeed the case. Consider assumptions of
Lemmas \ref{lem:boundedness_unstable} and
\ref{lem:boundedness_unstable_2}. In the proof of these lemmas
conditions (\ref{eq:stable}), (\ref{eq:g:Lipschitz:monotone:1})
played crucial role. These conditions provided estimates of the
upper bounds for $\dot{V}$, $-\dot{\lambda}$ as functions of
$V$ and $\lambda$. These upper bounds are used in conditions
(\ref{eq:positive:invariance:attracting}),
(\ref{eq:positive:invariance:attracting_2}) which are shown to
be incompatible with the assumption that the solutions passing
through the interior of $\Omega_a$ can cross the boundary
$\lambda=\psi(V(x))$. Thus one may expect that a dual result
might involve a set of conditions which is dual to
(\ref{eq:stable}), (\ref{eq:g:Lipschitz:monotone:1}). For this
purpose we consider the following substitutes for Assumptions
\ref{assume:stable},\ref{assume:unstable}:

\begin{assume}\label{assume:stable:2}  There exists  a continuous function $V:\Real^n\rightarrow \Real$, differentiable everywhere except, possibly, at the origin, and five functions of one variable,  $\underline{\alpha},\bar{\alpha}\in\mathcal{K}_{\infty}$, $\alpha: \ \Real_{\geq 0}\rightarrow\Real$, $\alpha(0)=0$, $\beta: \Real_{\geq 0}\rightarrow\Real_{\geq 0}$,
$\beta\in\mathcal{C}^{0}([0,\infty))$, $\varphi\in\mathcal{K}_0$,
such that  for every $(x,\lambda,t)\in{\mathcal{D}_{\Omega}}$ the following properties hold:
\begin{equation}\label{eq:stable:2}
\begin{split}
\underline{\alpha}(\|x\|)\leq V(x)\leq \bar{\alpha}(\|x\|), \ \
\frac{\pd V}{\pd x}f(x,\lambda,t)\geq  \alpha(V(x)) - \beta(V(x))\varphi(|\lambda|).
\end{split}
\end{equation}
\end{assume}
Assumption \ref{assume:stable:2} is similar to Assumption
\ref{assume:stable}, except for the sign of the last inequality
and also for the sign with which the term
$\beta(V)\varphi(|\lambda|)$ enters the right-hand side of
(\ref{eq:stable:2}). In essence, it states that  one can
provide a lower bound for the derivative of $V$ as a function
of $V$ and  $\lambda$. Let us  proceed with determining a
substitute for Assumption \ref{assume:unstable}.
\begin{assume}\label{assume:unstable:2}  There exist a function $\delta \in\mathcal{K}_0$ such that the following inequality holds for all $(x,\lambda,t)\in{\mathcal{D}_{\Omega}}$:
\begin{equation}\label{eq:g:Lipschitz:monotone:1:2}
\begin{split}
 g(x,\lambda,t)\leq - \delta(\|x\|).
\end{split}
\end{equation}
\end{assume}
Similarly to the case of Assumptions
\ref{assume:stable},\ref{assume:stable:2}, Assumption
\ref{assume:unstable:2} is almost a copy of Assumption
\ref{assume:unstable} in which the sign of the inequality is
reversed. The other difference is that we suppose that there is
an upper bound for $g(\cdot,\cdot,\cdot)$, and this bound  is
independent of $\lambda$. It is clear, however, that if such
dependence were to be in the form $g(x,\lambda,t)\leq
-\delta(\|x\|)-\xi(|\lambda|)$, $\delta,\xi\in\mathcal{K}_0$
then it can be reduced to the one stated in
(\ref{eq:g:Lipschitz:monotone:1:2}).

\begin{lem}[Solutions escaping a neighborhood of the origin]\label{lem:unboundedness_unstable}
Let system (\ref{eq:problem}) be given and satisfy Assumptions
\ref{assume:stable:2}, \ref{assume:unstable:2}. Suppose that
$\mathcal{D}$ contains the origin. In addition, let

\begin{itemize}
\item[(C5)] there exist  a function $\psi: \
    \psi\in\mathcal{K}_\infty
    \cap\mathcal{C}^{1}((0,\infty))$,
and an $a\in\Real_{>0}$ such that for some
$\varepsilon\in\Real_{>0}$ and all $V\in(0,a]$ the following  holds
\begin{equation}\label{eq:positive:invariance:repelling}
\begin{split}
\frac{\pd \psi(V)}{\pd V}&\left[\alpha(V)-\beta(V)\varphi(\psi(V)-\varepsilon)\right]+\delta\left(\overline{\alpha}^{-1}(V)\right)\geq 0;
\end{split}
\end{equation}
\item[(C6)] the set $\Omega_a$, where
\begin{equation}\label{eq:positive:invariance:domain:2}
\Omega_a=\{(x,\lambda) \ | \ x\in\overline{\mathcal{D}},  \
\psi(a)-\varepsilon\geq \lambda\geq \psi(V(x))-\varepsilon, \ V(x)\in[0,a]\}
\end{equation}
is contained in ${\mathcal{D}}\times (c_1,c_2)$.
\end{itemize}

Let $\Omega_a^{\ast}$ be the complement of $\Omega_a$ in
$\overline{\mathcal{D}}\times[c_1,\psi(a)-\varepsilon]$. Then
solutions of (\ref{eq:problem}) starting in $\Omega_a^\ast$ at $t=t_0$
cannot converge to the origin without  leaving the set
$\overline{\mathcal{D}}\times[c_1,\psi(a)-\varepsilon]$ at
least once.
\end{lem}

{\it Proof of Lemma \ref{lem:unboundedness_unstable}}.
{ Let
$\phi(\cdot;t_0,x_0,\lambda_0)$ $=$ $(x(\cdot;t_0,x_0,\lambda_0),
\lambda(\cdot;t_0,x_0,\lambda_0))$, $(x_0,\lambda_0)\in\Omega_a^\ast$ be a solution of (\ref{eq:problem}) converging to
the origin as in forward time   without leaving  the set
$\overline{\mathcal{D}}\times[c_1,\psi(a)-\varepsilon]$. This
implies that for an $r>0$ (sufficiently small) there is a
$t_1>t_0$ such that
$\lambda^2(t_1;t_0,x_0,\lambda_0)+\|x(t_1;t_0,x_0,\lambda_0)\|^2=r^2$,
and that  $\phi(t;t_0,x_0,\lambda_0) \in
\overline{\mathcal{D}}\times[c_1,\psi(a)-\varepsilon]$ for all $t\geq t_0$. Let
us pick
$r=(\overline{\alpha}^{-1}\circ\psi^{-1})(\varepsilon)/2$.
Notice that $\lambda(t_1;t_0,x_0,\lambda_0)$  is
positive since $\dot{\lambda}=g(x,\lambda,t)$ is non-positive
in $\mathcal{D}_\Omega$.
Thus the following estimate holds:
\[
\|x(t_1;t_0,x_0,\lambda_0)\|=\sqrt{r^2-\lambda(t_1;t_0,x_0,\lambda_0)} < (\overline{\alpha}^{-1}\circ\psi^{-1})(\varepsilon)<(\overline{\alpha}^{-1}\circ\psi^{-1})(\varepsilon+\lambda(t_1;t_0,x_0,\lambda_0)).
\]
This, however, implies that
$\psi(\overline{\alpha}(\|x(t_1;t_0,x_0,\lambda_0)\|))<
\lambda(t_1;t_0,x_0,\lambda_0)+\varepsilon$, and subsequently that
$\lambda(t_1;t_0,x_0,\lambda_0)>\psi(V(x(t_1;t_0,x_0,\lambda_0)))-\varepsilon$.

Introduce the function $p:[t_0,\infty)\rightarrow\Real$,
$p(t)=\psi(V(x(t;t_0,x_0,\lambda_0)))-\lambda(t;t_0,x_0,\lambda_0)-\varepsilon$.
The function $p(\cdot)$ is continuous and differentiable on $[t_0,t_1]$,   and  $\dot{p}(\cdot)$ is clearly bounded. Since $p(t_0)>0$ (condition C6) and $p(t_1)$<0
then, according to the intermediate value theorem, there must
be a point $t'\in(t_0,t_1)$ such that $p(t')=0$. Without loss
of generality we can suppose that $p(t)>0$ for all
$t\in[t_0,t')$. Let $t''$ be a point in $(t_0,t')$, then
\begin{equation}\label{eq:condition:contr}
p(t')-p(t'')=\int_{t''}^{t'}\dot{p}(s)ds=(t'-t'')\left(\frac{\pd \psi}{\pd V} \dot{V}(\tau)-\dot\lambda(\tau)\right)<0,
\end{equation}
where $\tau$ belongs to $[t'',t']$.   Recall that we assumed that
$\phi(t;t_0,x_0,\lambda_0)\in \overline{\mathcal{D}}\times[c_1,\psi(a)-\varepsilon]$ for all $t\geq t_0$, including
for $t\in[t'',t']$.
According to Assumptions
\ref{assume:stable:2} and \ref{assume:unstable:2}, the fact that
$\lambda(\tau;t_0,x_0,\lambda_0)\leq
\psi(V(x(\tau;t_0,x_0,\lambda_0)))-\varepsilon$, and condition C5
the following should hold:
\[
\begin{split}
&\frac{\pd \psi}{\pd V} \dot{V}(\tau)-\dot\lambda(\tau)\geq  \frac{\pd \psi}{\pd V}[\alpha(V(x(\tau;t_0,x_0,\lambda_0)))-\beta(V(x(\tau;t_0,x_0,\lambda_0)))\varphi(|\lambda(\tau;t_0,x_0,\lambda_0)|)]\\
&+\delta(\|x(\tau;t_0,x_0,\lambda_0)\|)\geq  \frac{\pd \psi}{\pd V}[\alpha(V(x(\tau;t_0,x_0,\lambda_0)))\\
&-\beta(V(x(\tau;t_0,x_0,\lambda_0)))\varphi(\psi(V(x(\tau;t_0,x_0,\lambda_0)))-\varepsilon)] + \delta(\overline{\alpha}^{-1}(V(x(\tau;t_0,x_0,\lambda))))\geq  0
\end{split}
\]
This, however, contradicts to (\ref{eq:condition:contr}). Hence
conditions of the lemma are incompatible with the assumption
that   the solution   $\phi(\cdot;t_0,x_0,\lambda_0)$ converges to the origin without
leaving the set
$\mathcal{D}\times[c_1,\psi(a)-\varepsilon]$ at least once.}
$\square$

When the function $\psi(\cdot)$ is not differentiable it is still
possible to provide conditions for specifying
solutions escaping a neighborhood of the origin at least once.
The conditions are dual to that of Lemma
\ref{lem:boundedness_unstable_2}, and we present them in the
lemma below.

\begin{lem}[Solutions escaping the origin 2]\label{lem:unboundedness_unstable_2}
Let system (\ref{eq:problem}) be given and satisfy Assumptions
\ref{assume:stable:2} and \ref{assume:unstable:2}. Suppose that
$\mathcal{D}$ contains the origin. In addition, let

\begin{itemize}
\item[(C7)] there exist  a function
    $\psi\in\mathcal{K}_\infty$ and an $a\in\Real_{>0}$
    such that for some $\varepsilon\in\Real_{>0}$ and all $V\in(0,a]$ the
    following holds
\begin{equation}\label{eq:positive:invariance:repelling:2}
\begin{split}
\psi(V)\left[\alpha(V)-\beta(V)\varphi(\psi(V)-\varepsilon)\right]+V\delta\left(\overline{\alpha}^{-1}(V)\right)\geq 0,
\end{split}
\end{equation}
and $\mathrm{hyp}(\psi_{[0,a]})$ is star-shaped with
respect to the origin.
\end{itemize}

\noindent Moreover, let condition (C6)   of Lemma
\ref{lem:unboundedness_unstable} holds,
and let  $\Omega_a^{\ast}$ be the complement of $\Omega_a$ in
$\overline{\mathcal{D}}\times[c_1,\psi(a)-\varepsilon]$.

Then solutions of (\ref{eq:problem}) with starting in
$\Omega_a^\ast$ at $t=t_0$ cannot converge to the origin without  leaving
the set
$\overline{\mathcal{D}}\times[c_1,\psi(a)-\varepsilon]$ at
least once.
\end{lem}

{\it Proof of Lemma \ref{lem:unboundedness_unstable_2}.} As before, let $\phi(\cdot;t_0,x_0,\lambda_0)$, $(x_0,\lambda_0)\in\Omega_a^\ast$ be a solution of (\ref{eq:problem})
converging to the origin as $t\rightarrow\infty$   without leaving the set
$\overline{\mathcal{D}}\times[c_1,\psi(a)-\varepsilon]$ at
least once. Invoking the same argument as in Lemma
\ref{lem:unboundedness_unstable} we can conclude that  there is
an interval $[t',t'']$ such that $\psi(V(x(t';t_0,x_0,\lambda_0)))-\lambda(t';t_0,x_0,\lambda_0)-\varepsilon=0$, and for all
$t\in[t',t'']$ we have that $\phi(t;t_0,x_0,\lambda_0)\in
\mathcal{D}\times \Lambda$ and
\begin{equation}\label{eq:hypoepi}
\begin{split}
\psi(V(x(t;t_0,x_0,\lambda_0)))-\lambda(t;t_0,x_0,\lambda_0)-\varepsilon&<0.
\end{split}
\end{equation}
It is clear that if the set $\mathrm{hyp}(\psi_{[0,a]})$ is
star-shaped with respect to the origin, then the set
$\mathrm{hyp}(\psi_{[0,a]}-\varepsilon)$ is star-shaped with
respect to the point  $(0,-\varepsilon)$.

Consider the function   $p: [t_0,\infty)\rightarrow\Real$,
$p(t)={V(x(t;t_0,x_0,\lambda_0))}/({\lambda(t;t_0,x_0,\lambda_0)+\varepsilon})$. 
The function $p(\cdot)$ is well defined for $t\in[t',t'']$. Notice
that, according to  (\ref{eq:hypoepi}), the point
$(V(x(t';t_0,x_0,\lambda_0)),\lambda(t';t_0,x_0,\lambda_0))\in
\mathrm{hyp}(\psi_{[0,a]}-\varepsilon)$, and
$(V(x(t'';t_0,x_0,\lambda_0)),\lambda(t'';t_0,x_0,\lambda_0))$ is in
the interior of $\mathrm{epi}(\psi_{[0,a]}-\varepsilon)$.
Therefore,  $p(t'')<p(t')$ for otherwise the point
$(V(x(t'';t_0,x_0,\lambda_0))$ , $\lambda(t'';t_0,x_0,\lambda_0))$
would be in $\mathrm{hyp}((\psi(V)-\varepsilon)_{[0,a]})$; this
would then contradict the condition
$\psi(V(t'';t_0,x_0,\lambda_0))-\lambda(t'';t_0,x_0,\lambda_0)-\varepsilon<0$.
On the other hand, there is a $\tau\in[t',t'']$ such that
\[
\begin{array}{ll}
p(t'')-p(t')=&(t''-t')[\dot{V}(x(\tau;t_0,x_0,\lambda_0))(\lambda(\tau;t_0,x_0,\lambda_0)+\varepsilon)\\
&-\dot{\lambda}(\tau;t_0,x_0,\lambda_0)V(x(\tau;t_0,x_0,\lambda_0))] {(\lambda(\tau;t_0,x_0,\lambda_0)+\varepsilon)^{-2}}.
\end{array}
\]
Hence, invoking Assumption \ref{assume:unstable:2} we arrive at
the following estimate:
\[
\begin{split}
&\dot{V}(x(\tau;t_0,x_0,\lambda_0))(\lambda(\tau;t_0,x_0,\lambda_0)+\varepsilon)-\dot{\lambda}(\tau;t_0,x_0,\lambda_0)V(x(\tau;x_0,\lambda_0))\geq\\
&\psi(V(x(\tau;t_0,x_0,\lambda_0)))\times\big[\alpha(V(x(\tau;t_0,x_0,\lambda_0)))-\beta(V(x(\tau;t_0,x_0,\lambda_0)))\times\\
&\times\varphi(\psi(V(x(\tau;t_0,x_0,\lambda_0)))-\varepsilon)\big]+V(x(\tau;t_0,x_0,\lambda_0))\delta\big(\overline{\alpha}^{-1}(V(x(\tau;t_0,x_0,\lambda_0))\big)\geq 0.
\end{split}
\]
This implies that $p(t'')-p(t')\geq 0$, and hence we have
reached a contradiction. $\square$

\section{Examples}\label{sec:Examples}

In this section we illustrate the theoretical results with examples. We begin with
the analysis of
systems (\ref{eq:simple_system:1})--(\ref{eq:simple_system:3})
presented in Example \ref{example:systems_toy} in Section \ref{sec:intro}.  Then we proceed to  the
phase synchronization problem introduced in Example
\ref{example:0}. Finally, in Section \ref{subsec:Applications}, we illustrate
how Lemma \ref{lem:boundedness_unstable_2} can be used to
approach an adaptive control problem for a class of systems nonlinear in the parameters.

\subsection{Forward invariant sets and basins of attraction for (\ref{eq:simple_system:1})--(\ref{eq:simple_system:3})}\label{sec:examples:academic}

\begin{figure}
\centering
\includegraphics[width=\textwidth]{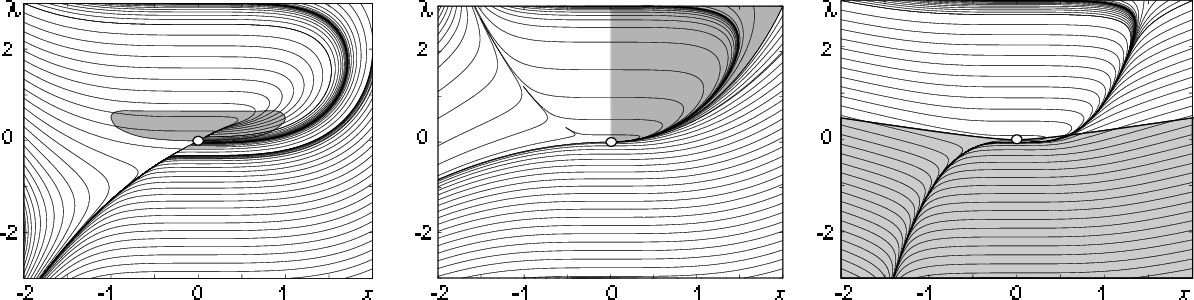}
\caption{{\it Left and middle panels:} forward invariant sets for systems
(\ref{eq:simple_system:1}), (\ref{eq:simple_system:2})  at
$\gamma=0.5$  (light grey areas in the plots).  The left panel
depicts the union of $\Omega_a(p)$, $p=(0,\infty)$, defined  by
(\ref{eq:examples:simple:1:domains}),  superimposed on the
phase plot of (\ref{eq:simple_system:1}). The middle panel shows
$\Omega_a(0.5)$, (\ref{eq:examples:simple:2:domains}), placed
over the phase plot of
(\ref{eq:simple_system:2}). {\it The right panel} depicts the domain (grey area)
 corresponding to solutions of (\ref{eq:simple_system:3}) not converging to the origin.}\label{fig:examples:toy}
\end{figure}

We start with {\it system (\ref{eq:simple_system:1})}:
\[
\begin{split}
\dot{x}&=-x + \lambda\\
\dot{\lambda}&=-\gamma |x|^3, \ \gamma\in\Real_{>0}.
\end{split}
\]
Let $\mathcal{D}=\Real$, $\Lambda=\Real$ and
notice that Assumption $\ref{assume:stable}$ is satisfied for
the first equation in (\ref{eq:simple_system:1}) with
$V(x)=x^2$, and $\alpha(V)=-2 V$, $\beta(V)=2\sqrt{V}$,
$\varphi(|\lambda|)=|\lambda|$. Assumption
\ref{assume:unstable} is fulfilled with $\xi(\cdot)\equiv 0$,
$\delta(|x|)=\gamma |x|^3$. Pick a candidate for the
function  $\psi: \ \psi(V)=p V$, $p\in\Real_{>0}$, and consider the
function $\mathcal{F}:\Real_{\geq 0}\rightarrow\Real$:
\begin{equation}\label{eq:B}
\mathcal{F}(V)=\frac{\pd \psi}{\pd V}(\alpha(V)+\beta(V)\varphi(\psi(V)))+\delta(\sqrt{V})=(-2 p + (2p^2+\gamma)\sqrt{V})V.
\end{equation}
The function is non-positive for
$\sqrt{V}\in(0,2p/(2p^2+\gamma)]$. Finally, note that  $\underline{\omega}(\psi(a))=\Real$ for all $a>0$. Hence, according to Lemma
\ref{lem:boundedness_unstable}, the sets
\begin{equation}\label{eq:examples:simple:1:domains}
\Omega_a(p)=\{(x,\lambda)\ | \ x\in\Real,\ \lambda\in\Real, \ p \left(\frac{2p}{2p^2+\gamma}\right)^2 \geq \lambda \geq p |x|^2, \ p\in\Real_{>0} \}
\end{equation}
are forward-invariant. Moreover,
$\lim_{t\rightarrow\infty}x(t;x_0,\lambda_0)=0$ for
$(x_0,\lambda_0)\in\Omega_a(p)$. Noticing that the right-hand side of (\ref{eq:simple_system:1}) is locally Lipschitz  and using Barbalat's lemma
we can conclude that
$\lim_{t\rightarrow\infty}\lambda(t;x_0,\lambda_0)=0$.  This
confirms that the origin is a weak attractor for
(\ref{eq:simple_system:1}). An estimate of the
forward invariant sets for (\ref{eq:simple_system:1}) is provided
in Fig. \ref{fig:examples:toy}. Note also that, since $\mathcal{D}=\Real$ and $\Lambda=\Real$, Corollary \ref{cor:boundedness_unstable} could be applied to this example too.

Let us now proceed with the analysis of {\it system
(\ref{eq:simple_system:2})}. Its dynamics is described as
\[
\begin{split}\dot{x}&=-x^2 + \lambda\\
\dot{\lambda}&=-\gamma |x|^3, \ \gamma\in\Real_{>0}.
\end{split}
\]
In this particular example we let $\mathcal{D}=\Real_{>0}$,
and $\Lambda=\Real$. One can see that  Assumptions
\ref{assume:stable} and \ref{assume:unstable} hold for this system
with $V=x^2$, $\alpha(V)=-2 V^{3/2}$, $\beta(V)=2\sqrt{V}$,
$\varphi(|\lambda|)=|\lambda|$, $\xi(\cdot)\equiv 0$, and
$\delta(|x|)=\gamma |x|^3$. Substituting these into
(\ref{eq:B})  and letting $\psi(V)=p V$ yields: $\mathcal{F}(V)=(-2p
+2p^2  +\gamma)V\sqrt{V}$.

It is therefore clear that $\mathcal{F}(V)\leq 0$ as long as $2p^2 - 2p
+\gamma\leq 0$. Let us suppose that the last inequality holds.
Hence condition (C1) of the lemma is satisfied. One can also
see that (C2) and (C3) hold too (with $\underline{\omega}(\psi(a))=\Real_{\geq 0}$ for all $a>0$). Thus Lemma
\ref{lem:boundedness_unstable} assures that the sets
\begin{equation}\label{eq:examples:simple:2:domains}
\Omega_a(p)=\{(x,\lambda)\ | \ x\in\Real,\ \lambda\in\Real, \ \lambda \geq p |x|^2, \ p: \  2p^2 - 2p +\gamma\leq 0, \ p\in\Real_{>0}\}
\end{equation}
are forward-invariant. Notice that $\gamma=0.5$ satisfies the
inequality $2p^2 - 2p +\gamma\leq 0$ for $p=0.5$. This is
consistent with the phase plots of (\ref{eq:simple_system:2})
in Example \ref{example:0} (see Fig.
\ref{fig:example:attractors}, $b$).

Finally, consider {\it(\ref{eq:simple_system:3})}:
\[
\left\{
\begin{array}{ll}\dot{x}&=-x^3 + \lambda\\
\dot{\lambda}&=-\gamma |x|^3.\end{array} \right.
\]
Let $\mathcal{D}=\Real$, and $\Lambda=\Real$.
The origin of this system is not an attractor. This can be shown
by employing the singular transformation
$(x,\lambda)\mapsto (\rho,\varphi)$: $x=\rho\cos(\varphi)$,
$\lambda=\rho^2\sin(\varphi)$ followed by the analysis of solutions
in  a vicinity of the point $(\rho,\varphi)=(0,0)$ in the $(\rho,\varphi)$ coordinates.

Let us see if we can determine domains corresponding to the
solutions of (\ref{eq:simple_system:3}) escaping some specified
neighborhoods of the origin.  For this purpose we will use
Lemma \ref{lem:unboundedness_unstable}. Let $V(x)=|x|^{q}$,
$q\in\Real_{>0}$, then Assumption \ref{assume:stable:2} holds
with $\alpha(\cdot),\beta(\cdot)$ and $\varphi(\cdot)$ defined as
\[
\alpha(V)=- q V^{\left(1+\frac{2}{q}\right)}, \  \beta(V)=q
V^{\left(1-\frac{1}{q}\right)}, \ \varphi(|\lambda|)=|\lambda|,
\]
and Assumption
\ref{assume:unstable:2} holds with $\xi(\cdot)\equiv 0$,
$\delta(|x|)=\gamma|x|^3$. Let $\psi(V)=pV^r$,
$r\in\Real_{>0}$, $p\in\Real_{>0}$, and consider
\[
\begin{split}
\mathcal{F}(V)=\frac{\pd \psi}{\pd V}(\alpha(V)-\beta(V)\varphi(\psi(V)-\varepsilon))+\delta(V^{\frac{1}{q}})
\geq -pqr V^{\left(r+\frac{2}{q}\right)}-p^2 qr V^{\left(2r-\frac{1}{q}\right)}+\gamma V^{\frac{3}{q}}.
\end{split}
\]
In order to apply Lemma \ref{lem:unboundedness_unstable} we
need to find $r,q,p\in\Real_{>0}$ and $a>0$ such that
$\mathcal{F}(V)\geq 0$ for $V\in(0,a]$. Letting $r=1$ and $q=2$ we arrive
at $\mathcal{F}(V)\geq -2p V^{2}-2 p^2 V^{\frac{3}{2}}+\gamma
V^{\frac{3}{2}}$. Hence $\mathcal{F}(V)\geq 0$ whenever $\sqrt V\leq
({\gamma}/{2p})-p$. Therefore, according to the lemma,
solutions of (\ref{eq:simple_system:3}) starting in the
complements $\Omega_a(p)^\ast$ of
\[
\Omega_a(p)=\{(x,\lambda) \ | \ x\in\Real, \ \left(\frac{\gamma}{2p}-p\right)^2-\varepsilon\geq\lambda\geq p V -\varepsilon, \ p\in(0,\sqrt{\gamma/2}) \}
\]
in
$\Real\times\left(-\infty,\left(\frac{\gamma}{2p}-p\right)^2-\varepsilon\right]$
cannot converge to the origin without leaving
$\Omega_a(p)^\ast$ at least once. The union of
$\Omega_a(p)^\ast$ over $p$ for $\gamma=0.5$ and $\varepsilon=0.01$ is
shown in Fig. \ref{fig:examples:toy}, right panel. Note that in
this particular case, since    $g(x,\lambda,t)=-\gamma |x|^3$, is non-positive    for all
$x,\lambda,t\in\Real$, solutions starting in this union do not
converge to the origin.


\subsection{Phase synchronization of neuronal oscillators}\label{subsec:example:phase_synch}

We now proceed with the analysis of a somewhat more
complicated system: a network of coupled neuronal oscillators
(or cells). Interaction between individual elements is allowed
to be heterogeneous: the cells are able to interact with
immediate neighbors via gap-junctions (intercellular
connections enabling direct flows of ions from one to another),
or can transmit pulses to other cells through synaptic
connections.

\begin{figure}
\centering
\includegraphics[width=0.8\textwidth]{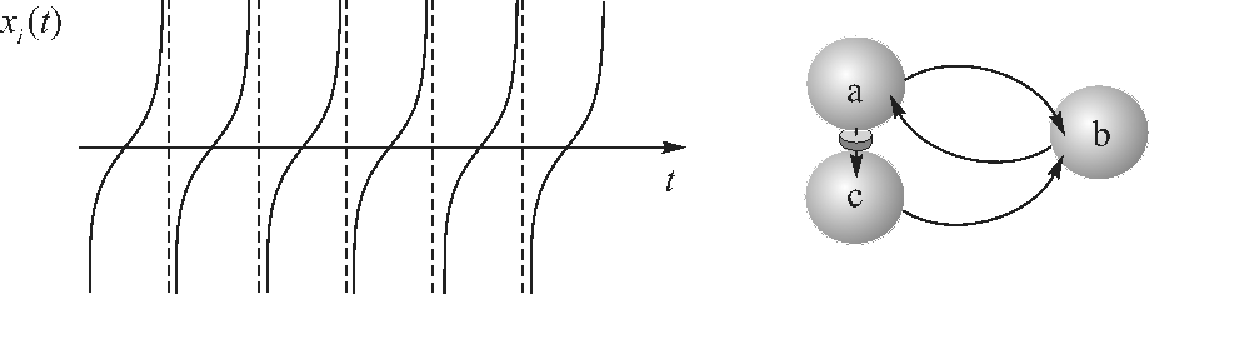}
\caption{Networks of phase oscillators. The left panel shows a diagram depicting solutions $x_i(t)=-\cot(t+\phi_i)$, $\phi_i\in\Real$, $i\in\{a,b,c\}$, $t\neq\pi k- \phi_i$, $k\in\mathbb{Z}$, produced by the $i$-th equation in the uncoupled system (\ref{eq:cell:1})--(\ref{eq:cell:3}) (i.e. at $\varepsilon=0$, $\varepsilon_1=0$) as as function of time.
The right panel  depicts a diagram of the network described by (\ref{eq:cell:1})--(\ref{eq:cell:3}). Long arrows indicate long-range connections via pulse coupling, and the short arrow indicates connection between cells $a$ and $c$ via unidirectional gap junction.}\label{fig:example:phase}
\end{figure}

In the absence of coupling each oscillator generates spikes
with a given frequency. The oscillators may not be identical, yet the frequency of oscillations is
supposed to be the same for every element in the network. For
illustrative purposes we will assume that oscillations in the
cells occur through the saddle-node-on-invariant-cycle
bifurcation. Such mechanism is inherent to a wide range of
models, including the canonical Hodgkin-Huxley equations describing
potential generation in neural membranes via potassium-sodium
gates cf. \cite{Hodgkin_Huxley}, \cite{Izhikevich:2007}.
Oscillators of this kind have the following normal form
\cite{Izhikevich:2007}:
\begin{equation}\label{eq:normal_form}
\dot{x}_i=1 + x_i^2 + \varepsilon_i u(t),
\end{equation}
where $x_i$ represents the neuron's membrane potential,
$\varepsilon_i\in\Real_{\geq 0}$ is an input gain, and
$u:\Real\rightarrow\Real$ is a function that models the inputs
(couplings) applied to the cell. In this model, when $x_i(\cdot;t_0,x_0)$
escapes to infinity a spike is produced and the initial
conditions are reset to $-\infty$. The process repeats
infinitely many times giving rise to periodic spikes with infinite
amplitude (see Fig. \ref{fig:example:phase}).



In this particular example we consider a network consisting of
three coupled neuronal oscillators, (\ref{eq:normal_form}), of
which the topological structure is shown in Fig.
\ref{fig:example:phase}. Neurons $a$ and $b$, interact via
bi-directional pulse coupling, neuron $a$ interacts with neuron
$c$ via a gap-junction, and $c$ interacts with $a$ synaptically
(pulse coupled). Dynamics of such a network can be described as
follows \cite{Izhikevich:2007}:
\begin{align}
    \dot x_a &= 1+x_a^2 + \varepsilon[\delta (t-t_b)+\delta(t-t_c)], \label{eq:cell:1} \\
    \dot x_b &= 1+x_b^2 + \varepsilon\delta (t-t_a),     \label{eq:cell:2} \\
    \dot x_c &= 1+x_c^2 + \varepsilon_1(x_a-x_c),     \label{eq:cell:3}
\end{align}
where $\delta(\cdot)$ is the Dirac's delta function. When
$\varepsilon,\varepsilon_1$ are sufficiently small  phases
of each oscillator can be defined as follows:
$\varphi_i=t+\phi_i$,  $i\in\{a,b,c\}$, where $\phi_i$ are
``slow'' phase fluctuations. Despite the fact that the frequencies of  the individual
oscillators are kept identical, phase fluctuations, $\phi_i$,
may vary. The question is: whether these fluctuations in the
coupled system, as functions of $t$ converge to some known values as
$t\rightarrow\infty$.

In order to answer to this question consider the dynamics of slow
phase fluctuations in system
(\ref{eq:cell:1})--(\ref{eq:cell:3}) (cf.
\cite{Izhikevich:2007}):
\begin{align}
    \dot \phi_a &= \varepsilon/\pi[\sin^2(\phi_b-\phi_a) + \sin^2(\phi_c-\phi_a)], \nonumber  \\
    \dot \phi_b &= \varepsilon/\pi\sin^2(\phi_a-\phi_b), \nonumber  \\
    \dot \phi_c &= \varepsilon_1/2\sin(2(\phi_a-\phi_c)). \nonumber
\end{align}
Denoting $\lambda:= \phi_b-\phi_a$ and $x:= \phi_a-\phi_c$ this
results in:
\begin{align}
       \dot x &= \tilde \varepsilon[\sin^2(\lambda) + \sin^2(x)] - \varepsilon_1/2\sin(2x), \ \tilde \varepsilon:= \frac{\varepsilon}{\pi},
 \label{eq:phase:2}\\
  \dot \lambda &=  - \tilde \varepsilon \sin^2(x),  \label{eq:phase:1}
\end{align}
Note that equations  (\ref{eq:phase:2}), (\ref{eq:phase:1})
locally resemble system (\ref{eq:problem}). Indeed, the origin
of (\ref{eq:phase:2}) at $\lambda=0$ is  locally asymptotically
stable in the region $\{x\in(-\pi/2,\pi/2) \vert
\frac{\varepsilon_1}{\tilde \varepsilon}> \tan(|x|)\}$, and the
right-hand side of (\ref{eq:phase:1}) is a non-negative
function of $x$, $\lambda$.

Consider $V(x) = \tan^2(x)$, then for all $|x|<M$, $M=\tan^{-1}\left|
\frac{\varepsilon_1}{\tilde \varepsilon}
\right|$,  $\lambda\in\Real$  the
following holds: $\dot V \leq \alpha(V) + \beta(V)\varphi(|\lambda|)$,
where
\[
\begin{split}
    \alpha(V) &= -2 (\varepsilon_1-\tilde \varepsilon V^{1/2})V,  \\
    \beta(V) &= 2 \tilde \varepsilon V^{1/2}(1+V),
\end{split}
\ \ \
\begin{split}
    \varphi(|\lambda|) = \left\{\begin{array}{l}\sin^2(|\lambda|), \ |\lambda|\in [0,\pi/2)\\
    1, \ |\lambda|\geq \pi/2. \end{array}  \right.
    \end{split}
\]
Thus letting $\mathcal{D}=(-\mu,\mu)$, $0<\mu<M$,
$\Lambda=\Real$, and $\underline{\alpha}(|x|)
=\overline{\alpha}(|x|)= \tan^{2}(|x|)$ we can conclude that
Assumption \ref{assume:stable} holds for  (\ref{eq:phase:2}),
(\ref{eq:phase:1}). With regards to Assumption
\ref{assume:unstable}, it is satisfied with  $\delta(|x|) =
\tilde \varepsilon \varphi(|x|)$, $\xi(|\lambda|)\equiv 0$.

\begin{figure}
\centering
\includegraphics[width=240pt]{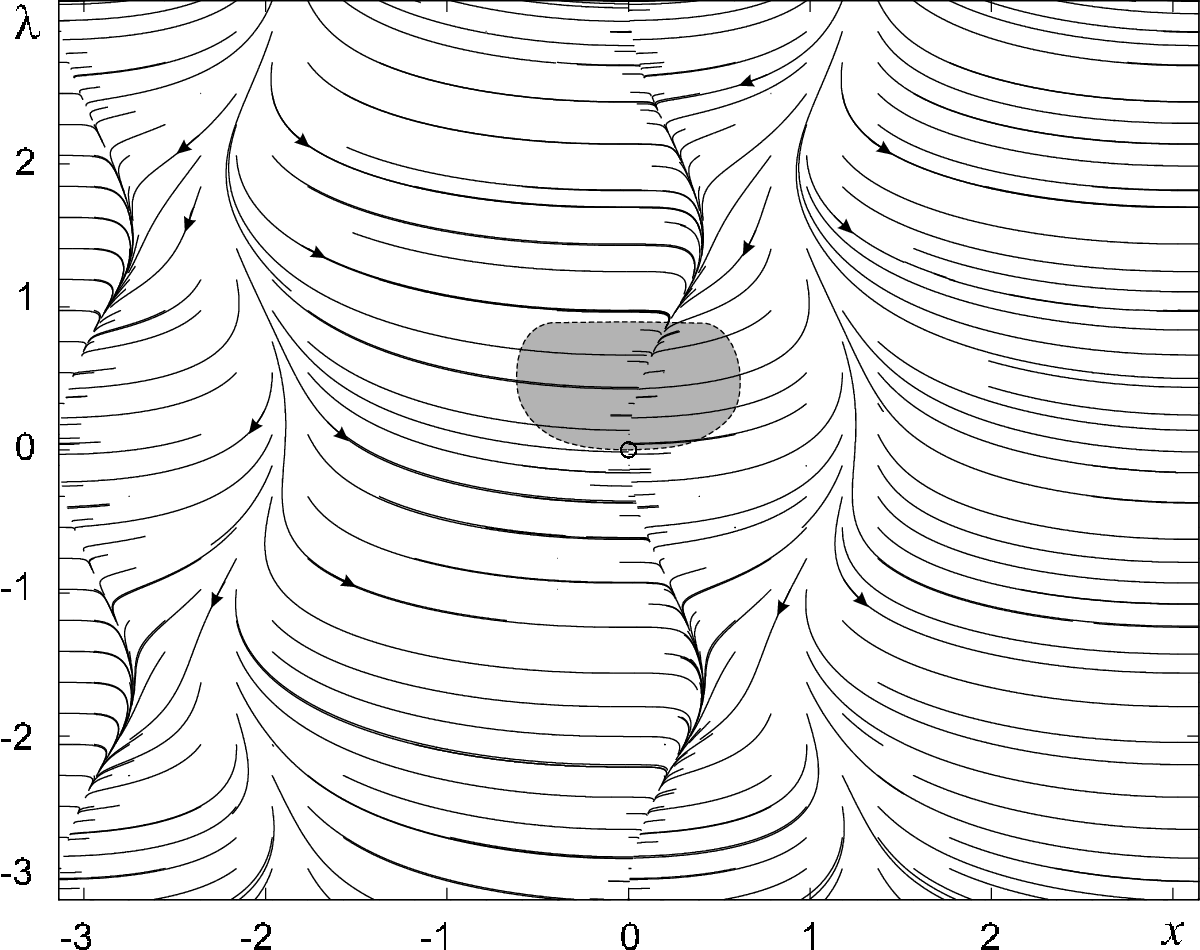}
\caption{Estimate of the forward invariant set for system (\ref{eq:phase:2}), (\ref{eq:phase:1})
at  $\varepsilon_1=\varepsilon=0.1$. The estimate (grey shaded area) is the union $\cup_{p\in(0,10]}\Omega_a(p)$  of the sets $\Omega_a(p)$ defined in (\ref{eq:phase:3:domains}). }\label{fig:phase_synch_example}
\end{figure}

Suppose that we wish to apply Lemma
\ref{lem:boundedness_unstable}. In this case we need to find a
strictly monotone function $\psi\in\mathcal{K}$ such that the
function $\mathcal{F}(\cdot)$ defined as

\begin{align}
\mathcal{F}(V)=\frac{\partial \psi}{\partial V}(\alpha(V)+\beta(V)\varphi(\psi(V)))+\delta(\underline{\alpha}^{-1}(V)), \label{eq:ineq}
\end{align}
is non-positive in $(0,a]$ for some $a\in\Real_{>0}$. Let
$\psi(\cdot)$ be such that its restriction on the set
$\{V\in\Real_{\geq 0} \ | \ V\leq \underline{\alpha}(\mu)\}$ is
as follows: $\psi(V) = p\tan^{-1}(V^{\tfrac{1}{2}}), \
p\in\Real_{>0}$. Then for all $V\in(0,\underline{\alpha}(\mu)]$
the function $\mathcal{F}(\cdot)$ in (\ref{eq:ineq}) can be estimated
from above as follows:
\[
\begin{split}
\mathcal{F}(V)&=p\frac{-2 (\varepsilon_1-\tilde \varepsilon V^{\frac{1}{2}})V+2 \tilde \varepsilon V^{\frac{1}{2}}(1+V) \sin^2(p\tan^{-1}(V^{\frac{1}{2}}))}{2V^{\frac{1}{2}}(1+V)}+\tilde\varepsilon\sin^2(\tan^{-1}(V^{\frac{1}{2}}))\\
&\leq -\frac{p\varepsilon_1 V^{\frac{1}{2}}}{(1+V)}+\tilde{\varepsilon}(p^3+p+1) V.
\end{split}
\]
Hence $\mathcal{F}(V)\leq 0$ whenever $V\leq \tan^2(\mu)$ and $-p\varepsilon_1+\tilde{\varepsilon}(p^3+p+1) V^{\frac{1}{2}}(1+V)\leq 0$. Thus solving the latter inequality for $V^{\frac{1}{2}}$ and choosing $V\leq a(p)$,  where
\[
\begin{split}
a(p)=\min\left\{\tan^2(\mu), \left(r(p)-\frac{1}{3r(p)}\right)^2 \right\}, \ r(p)=\left(\frac{p\epsilon_1+\sqrt{(p\epsilon_1)^2+\frac{4}{27}\tilde\varepsilon^2(p^3+p+1)^2}}{2\tilde\varepsilon(p^3+p+1)}\right)^{\frac{1}{3}}, \end{split}
\]
we can ensure that condition (C1) of the lemma is satisfied.
Given that the ball $\{ x | \ x\in\Real^n, \ |x|\leq \tan^{-1}(a)\}$ is contained in
$\mathcal{D}$, (condition (C2)), and that
${ \underline{\omega}(\psi(a))}=\Real$, Lemma \ref{lem:boundedness_unstable}
implies that
\begin{equation}\label{eq:phase:3:domains}
    \Omega_a(p) = \left\{(x,\lambda)\vert \  p\tan^{-1}(a^{\frac{1}{2}}(p)) \geq \lambda \geq p |x| \right\}
\end{equation}
is forward-invariant. The union of $\Omega_a(p)$, $p\in(0,10]$ as well as the phase plot of (\ref{eq:phase:2}) for $\varepsilon_1=0.1$, $\varepsilon=0.1$ are shown in Fig. \ref{fig:phase_synch_example}.

Thus we can conclude that if the values of relative phases $\phi_{a,0}$, $\phi_{b,0}$, and $\phi_{c,0}$ are chosen so that $(\phi_{b,0}-\phi_{a,0},\phi_{a,0}-\phi_{c,0})\in \cup_{p> 0} \Omega_a(p)$, where  $\Omega_a(p)$ is defined in (\ref{eq:phase:3:domains}), then the relative phases $\phi_b(\cdot;t_0,\phi_{a,0},\phi_{b,0},\phi_{c,0})-\phi_a(\cdot;t_0,\phi_{a,0},\phi_{b,0},\phi_{c,0})$ and $\phi_a(\cdot;t_0,\phi_{a,0},\phi_{b,0},\phi_{c,0})-\phi_c(\cdot;t_0,\phi_{a,0},\phi_{b,0},\phi_{c,0})$  are bounded on $[t_0,\infty)$. Moreover $\lim_{t\rightarrow\infty} \phi_a(t)-\phi_c(t)=0$, and consequently, $\lim_{t\rightarrow\infty} \phi_b(t)-\phi_a(t)=0$. Hence relative phases in the system converge to zero as $t\rightarrow \infty$ provided that $(\phi_{b,0}-\phi_{a,0},\phi_{a,0}-\phi_{c,0})\in \cup_{p> 0} \Omega_a(p)$. This corresponds to in-phase synchronization of the solutions of (\ref{eq:cell:1})--(\ref{eq:cell:3}).  Looking at the phase plot in Fig. \ref{fig:phase_synch_example} one can observe that the actual domain of initial conditions corresponding to in-phase synchronization is substantially larger than the one obtained analytically (grey area). This may be viewed as a consequence of choosing functions $V(\cdot)$, $\psi(\cdot)$, $\alpha(\cdot)$ and $\beta(\cdot)$ so that the derivations of $\Omega_a(p)$  are kept simple.

\subsection{Adaptive control for  nonlinearly parameterized systems}\label{subsec:Applications}

As a yet another and final illustration consider the following
system:
\begin{equation}\label{eq:system_ac}
\dot{x}=f(x,\theta,t,u(t)),
\end{equation}
where the function
$f:\Real^n\times\Real^m\times\Real\times\Real^q\rightarrow\Real^n$
is continuous and locally Lipschitz,
$\theta\in \Theta$,
$\Theta\subset\Real^m$ is the vector of unknown
parameters, and $u:\Real\rightarrow\Real$ is a control input. With regards to the
input $u$, we suppose that it is a continuous function of $x$
and a parameter $\mu\in\Real$. Furthermore, we suppose that the
  for any $\theta\subset\Theta$ there is a stabilizing feedback $u(t)=\tilde{u}(x(t),\mu)$  that is locally Lipschitz in $\mu$. In particular,   we suppose that the following assumption holds (cf \cite{CDC:Pomet:1992}):

\begin{assume}\label{assume:control} Let ${\mathcal{M}}=[\mu_1,\mu_2]$, $\mu_2>\mu_1$ be an interval in $\Real$, and let $V:\Real^n\rightarrow\Real$ be a positive-definite and differentiable function:
\[
\underline{\alpha}(\|x\|)\leq V(x)\leq \overline{\alpha}(\|x\|), \ \underline{\alpha}, \ \overline{\alpha}\in\mathcal{K}_\infty.
\]
The function $f(\cdot,\cdot,\cdot,\cdot)$ in  (\ref{eq:system_ac}) is such that:

\begin{itemize}
 \item[ 1)] there are known continuous and locally Lipschitz functions
     $\tilde{u}:\Real^n\times\Real\rightarrow\Real^q$ and  $\alpha\in\mathcal{K}_\infty$   such that for
     any $\theta\in \Theta$ there is a
     $\mu\in \mathcal{M}$:
\begin{equation}\label{eq:asymptotic_stability}
\frac{\pd V}{\pd x}f(x,\theta,t,\tilde{u}(x,\mu))\leq - \alpha(V);
\end{equation}
\item[2) ] there is a $\beta\in\mathcal{K}$ such that for
    every $\mu,\mu'\in \mathcal{M}$ and
    $\theta\in \Theta$ the following holds:
\[
\frac{\pd V}{\pd x}(f(x,\theta,t,\tilde u(x,\mu))-f(x,\theta,t,\tilde u(x,\mu')))\leq \beta(V)|\mu-\mu'|.
\]
\end{itemize}
\end{assume}
Checking that the assumption holds for a given particular
system may present certain technical difficulties. Detailed
discussion of this issue is outside of the scope of this work.

We are interested in finding a function
$\hat{\mu}:\Real\rightarrow\Real$ such that solutions of
$\dot{x}=f(x,\theta,t,\tilde u(x,\hat{\mu}(t)))$ converge to the
origin as $t\rightarrow\infty$ for all
$\theta\in \Theta$. Let  $\hat{\mu}$:
\begin{equation}\label{eq:adaptive_law}
\begin{split}
\hat{\mu}=q(h), \ \dot{h}=-\gamma(\|x\|), \ \gamma\in\mathcal{K},
\end{split}
\end{equation}
where the function $q:\Real\rightarrow\Real$ is  $H$-periodic, that is $q(h+H)=q(h)$ for all $h\in\Real$, and
Lipschitz with constant $\ell>0$,  $\max_{h\in\Real} q(h)=\mu_2$,
$\min_{h\in\Real} q(h)=\mu_1$.

Now we are ready to formulate the following result:
\begin{cor}\label{thm:adaptive} Consider system (\ref{eq:system_ac}) and let it satisfy Assumption \ref{assume:control}. Consider  (\ref{eq:system_ac}), (\ref{eq:adaptive_law}):
\begin{equation}\label{eq:closed_loop_adaptive}
\begin{split}
\dot{x}=f(x,\theta,t,\tilde u(x,\hat{\mu})), \ \
\hat{\mu}=q(h), \ \
\dot{h}=-\gamma(\|x\|),
\end{split}
\end{equation}
and suppose that the following conditions hold:
\begin{itemize}
\item[(C8) ] there is a $V^\ast>0$ such that
    $\alpha(V)>\beta(V)(\mu_2-\mu_1) \ \mbox{ for all }\ V\geq
    V^\ast$;
\item[(C9) ] $\mathrm{epi}(\alpha_{[0,V^\ast]})$ is star-shaped w.r.t. the origin, and the function $\gamma(\cdot)$ in
    (\ref{eq:closed_loop_adaptive}) satisfies:
\begin{equation}\label{eq:gamma_choice}
\gamma\in\mathcal{K}, \ \gamma(\underline{\alpha}^{-1}(V))<\frac{\alpha(V)^2}{4\beta(V^\ast)\ell V}.
\end{equation}
\end{itemize}
Then, solutions of (\ref{eq:closed_loop_adaptive})  are
bounded, and furthermore
\begin{equation}\label{eq:theorem_conclusion}
\lim_{t\rightarrow\infty}x(t;t_0,x_0,h_0)=0, \ \mbox{and \ there \ is \ a} \ \mu'\in\mathcal{M}: \  \lim_{t\rightarrow\infty}\hat{\mu}(q(t;t_0,x_0,h_0))=\mu'.
\end{equation}
\end{cor}
{\it Proof of Corollary \ref{thm:adaptive}}. Let $x_0\in\Real^n$, $h_0\in\Real$. We begin with
showing that solutions of the interconnection are defined for
all $t\geq t_0$.  Given that the right-hand side of
(\ref{eq:closed_loop_adaptive}) is locally Lipschitz, the solution satisfying the initial condition $x(t_0)=x_0$, $h(t_0)=h_0$ is defined at least locally in a neighborhood of $(x_0,h_0,t_0)$.
Let $[t_0,t_1)$, $t_1<\infty$ be the maximal interval of the
solution's definition for $t\geq t_0$. This implies that for  any $M\in\Real_{\geq 0}$ there exists a
$t'\in[t_0,t_1)$ such that
$\max\{|h(t';t_0,x_0,h_0)|,\|x(t';t_0,x_0,h_0)\|\} \geq M$. According to Assumption
\ref{assume:control}, there is a $\mu\in \mathcal{M}$ such
that:
\[
\dot{V}\leq \frac{\pd V}{\pd x}f(x,\theta,t,\tilde u(x,\mu))+ \frac{\pd V}{\pd x}(f(x,\theta,t,\tilde u(x,\hat{\mu})-f(x,\theta,t,\tilde u(x,\mu))\leq -\alpha(V)+\beta(V)|\mu-\hat{\mu}|.
\]
Taking condition C8 of the theorem into account and using the
comparison lemma, we conclude that the function $V(x(\cdot;t_0,x_0,h_0))$ is
bounded from above by $\max\{V^\ast,V(x_0)\}$ on  $[t_0,t']$. This, however, contradicts
 $V(x(t';t_0,x_0,h_0))\geq
\underline{\alpha}(\|x(t';t_0,x_0,h_0)\|)\geq
\underline{\alpha}(M)$  because $M$  can be chosen so that $\underline{\alpha}(M) > \max\{V^\ast,V(x_0)\}$ . Therefore $x(\cdot;t_0,x_0,h_0)$ must be
bounded on $[t_0,t_1)$. On the other hand, according to
(\ref{eq:gamma_choice}), the following holds:
\[
|h(t';t_0,x_0,h_0)|\leq |h_0| + (t'-t_0)\max_{V\in[0,\max\{V(x_0),V^\ast\}]}\left\{\frac{\alpha(V)^2}{4\beta(V^\ast) \ell V} \right\}.
\]
Thus $|h(\cdot;t_0,x_0,h_0)|$ is bounded on $[t_0,t_1)$.
Hence solutions of (\ref{eq:closed_loop_adaptive}) are defined
for all $t \geq t_0$.  Furthermore, as follows from C8, the
function $V(x(\cdot;t_0,x_0,h_0))$  and, consequently,  $x(\cdot;t_0,x_0,h_0)$
are bounded on $[t_0,\infty)$. Moreover, there is a
$t^\ast\geq t_0$, independent on $h_0$, such that
$V(x(t;t_0,x_0,h_0))\leq V^\ast$ for all $t \geq t^\ast$.

Let us show that  $h(\cdot;t_0,x_0,h_0)$ is bounded on
$[t_0,\infty)$ as well. In order to do so we will invoke Lemma
\ref{lem:boundedness_unstable_2}. According to Assumption
\ref{assume:control} there is $\mu\in \mathcal{M}$ such
that: $(\pd V \ / \ \pd x) f(x,\theta,t,\tilde u(x,\mu))\leq
-\alpha(V)$. Given that the function $q(\cdot)$ in the definition of
$\hat{\mu}$ is continuous and periodic, it follows from the
intermediate value theorem that there is a $h_{\mu}\in(0,H)$:
$q(h_{\mu}+n H)=\mu$, for all $n\in\Numbers$. Let us pick an $n$ such
that
\begin{equation}\label{eq:searching_variable:1}
\eta(n)=h_\mu+n H, \ \eta(n)<h_0, \ h(t^\ast,x_0,h_0)-\eta(n)\geq \frac{\alpha(V^\ast)}{2\beta(V^\ast)\ell},
\end{equation}
and define $\lambda(t)=h(t;t_0,x_0,h_0)-\eta(n)$.  It is clear that
for $t\geq t^\ast$ the following holds along the solutions of
(\ref{eq:closed_loop_adaptive}):
\[
\begin{split}
\dot{V}\leq \frac{\pd V}{\pd x}f(x,\theta,t,\tilde u(x,\hat{\mu}))\leq -\alpha(V)+\beta(V^\ast)\ell|\lambda(t)|, \ \
-\gamma(\underline{\alpha}^{-1}(V))\leq\dot{\lambda}(t)\leq 0.
\end{split}
\]

Consider the following function $\psi\in\mathcal{K}_\infty$:
$\psi(V)={\alpha(V)}/({2\beta(V^\ast)\ell})$. According to Lemma
\ref{lem:boundedness_unstable_2} and Remark \ref{rem:simple_generalization} the set
\[
\Omega_a=\{(x,\lambda) \ | \ x\in\Real^n, \ \lambda\in\Real, \ \psi(a)\geq \lambda\geq\psi(V(x)), \ V(x)\in[0,a]\}
\]
is forward invariant for $t\geq t^\ast$ provided that
$\psi(V)(-\alpha(V)+\beta(V^\ast)\ell
\psi(V))+V\gamma(\underline{\alpha}^{-1}(V))\leq 0$  for all $
V\in[0,a]$. It is therefore clear that the choice:
$\gamma(\underline{\alpha}^{-1}(V))\leq
{\psi(V)\alpha(V)}/({2V})= {\alpha(V)^2}/({4\beta(V^\ast)V \ell})$
ensures that the set $\Omega_a$ is forward invariant. Hence
every solution of (\ref{eq:closed_loop_adaptive}) satisfying
$h(t^\ast;t_0,x_0,h_0)=\lambda(t^\ast)+\eta(n),  \ \lambda(t^\ast)=\frac{\alpha(V^\ast)}{2\beta(V^\ast)\ell}$, and $ x(t^\ast;t_0,x_0,h_0): \ V(x(t^\ast,x_0,h_0))\leq V^\ast$
should necessarily be bounded on $[t^\ast,\infty)$. This, in view of the choice of
$n$ in the definition of the variable $\lambda(t)$,
(\ref{eq:searching_variable:1}), assures that
$\lambda(\cdot)$ is bounded on $[t^\ast,\infty)$. Hence
$h(\cdot;t_0,x_0,h_0)$ is bounded on $[t_0,\infty)$. This implies
$\lim_{t\rightarrow\infty}\int_{t_0}^{t}\gamma(\|x(\tau;t_0,x_0,h_0)\|)d\tau
=h' < \infty$, where the function
$\gamma(\|x(\cdot;t_0,x_0,h_0)\|)$ is
uniformly continuous on $[t_0,\infty)$. Thus,  invoking Barbalat's lemma we can
conclude that (\ref{eq:theorem_conclusion}) holds. $\square$

Note that despite (\ref{eq:theorem_conclusion}) holds,  the values of the internal state $h$ of the controller can be large, depending on the initial conditions $x_0,h_0$. This is a well-known drawback of schemes of this type \cite{Dyn_Con:Ilchman:97,CDC:Pomet:1992}. One can remove this limitation by replacing (\ref{eq:adaptive_law}) with  $\hat{\mu}(h)=\mu_1+0.5(\mu_2-\mu_1)(1+h), \ \dot{h}=-\gamma(\|x\|) \frac{2\ell}{\mu_2-\mu_1} z , \ \dot{z}=\gamma(\|x\|)\frac{2\ell}{\mu_2-\mu_1} h, \ h(t_0)=h_0, \ z(t_0)=z_0, \ \gamma\in\mathcal{K}$ and restricting the initial condition $h_0$, $z_0$ to $h_0^2+z_0^2=1$ (see \cite{SIAM_non_uniform_attractivity}). Another issue is that the time needed for $x(\cdot;t_0,x_0,h_0)$ to converge to a given neighborhood of the origin may be large too. Derivation of a-priori estimates of the convergence times requires further analysis and, possibly, additional constraints.

\section{Conclusions}\label{sec:Conclusion}

In this manuscript we presented results for  finding
forward invariant sets and assessing convergence of
solutions in dynamical systems with unstable equilibria. In
particular, we focused on systems in which stable motions in
higher dimensional subspace of the system state space are
coupled with unstable motions in a lower dimensional subspace.
Such systems, as has been illustrated with examples, occur in a
relevant range of problems including adaptive control in the
presence of general nonlinear parametrization of uncertainty,
and phase synchronization in networks of coupled oscillators.

Motivated by some limitations of earlier analysis techniques proposed for this class of systems  (cf. \cite{SIAM_non_uniform_attractivity}, \cite{Dyn_Con:Ilchman:97}), we aimed at developing a more versatile alternative.  The alternative should apply to systems with unstable attractors while, at the same time, retain convenience and simplicity of  conventional Lyapunov-functions based analysis.

The method proposed and discussed in the article allows to produce simple algebraic tests for finding areas of forward invariance for systems with Lyapunov-unstable invariant sets. Moreover, the method can be applied to checking whether given equilibrium is an attractor, albeit not necessarily stable. Estimates of the attractor basins  are also supplied. In addition to convergence and boundedness criteria, geometric intuition  behind our results allows one to approach a dual problem: the one of estimating domains of initial conditions corresponding to solutions escaping given neighborhood of the origin at least once. These latter results are relevant in the context of determining relaxation times in nonlinear dynamical systems \cite{Gorban:2004}.

In spite of advantages of the method we note that there are limitations too. In particular, we require that the function $V(\cdot)$ characterizing dynamics of the stable part in (\ref{eq:problem}) vanishes only at a single point, $x=0$. This prevents explicit applications of the results to systems in which solutions of the stable part (in absence of coupling with the unstable part) are attracted to an orbit or a set which is not a single point.  We hope, however, to be able to address this and other issues in future publications.

\bibliographystyle{plain}
\bibliography{Unstable_convergence_lemma}

\section{Appendix. Star-shaped sets and envelopes}\label{sec:Appendix}

Let $V$ be a real vector space. In what follows symbol $[x,y]$
will denote a segment connecting two vectors $x,y \in V$:
$[x,y]=\{\gamma x + (1-\gamma) y \mid \gamma \in [0,1] \}$.

\begin{dfn}
A set $S \subset V$ is { \it star shaped} with respect to a
point $x \in S$, if for any $y\in S$ the   segment $[x,y]$  also
belongs to $S$: $[x,y] \subset S$.
\end{dfn}

\noindent The following properties hold for star shaped sets in
$V$:

$\bullet$ A set is convex iff it is star shaped  with respect
    to its every point.

$\bullet$ Let  $x \in V$ and $W$ be a family of sets star shaped
    with respect to $x$. Then both the intersection
    $\cap_{S\in W} S$  and the union
    $\cup_{S\in W} S$  are star shaped with
    respect to $x$.

$\bullet$ Let $E$ be a real vector space,  $A:V \to E$ be a
    linear map, and $S \subset V$ be a star shaped set with
    respect to a point $x \in S$. Then the image $A(S)$ of
    $A(\cdot)$ is star shaped with respect to $A(x)$.

\begin{dfn}
For any set $D \subset V$ and a point $x \in D$ the star shaped
envelope of $D$ with respect to $x$, ${\rm star}_x(D)$, is the
minimal star shaped set with respect to $x$ which includes $D$.
That is: every shaped set with respect to $x$ including $D$
must include $\mathrm{star}_x(D)$.
\end{dfn}

Star shaped envelopes exist and could be defined in two
alternative ways. Namely, ``from above'' (as an intersection):
\begin{equation}\label{a:star_shaped_envelop_set_above}
\begin{array}{c}
{\rm star}_x(D)=\bigcap_{S \in W_x(D)} S,
\end{array}
\end{equation}
where $W_x(D)$ is a family of sets that are star shaped with
respect to $x$ and include $D$, and  ``from below'' (as a union
of segments):
\begin{equation}\label{defstarenvset}
\begin{array}{c}
S=\bigcup_{y\in D}[x,y].
\end{array}
\end{equation}
Notice that deriving a star shaped envelope of an analytically
defined set computationally is a much easier procedure than
that of deriving a convex envelop of the same set.

Let us remind that epigraph of a real valued function
$f:S\rightarrow\Real$  is a subset of $S \times \Real$ that
consists of all points lying on or above its graph: ${\rm epi}
(f)=\{(x,\gamma) \mid x \in S, \, \gamma \geq f(x)\}.$

\begin{dfn}
  The function $f:S\rightarrow\Real$ is {\it star shaped} with respect to $x \in S$, if ${\rm epi} (f)$ is star shaped
set with respect to $(x,f(x))$. The function is {\it convex} if
it is star shaped with respect to every $x\in S$.
\end{dfn}
If  the function $f:S\rightarrow\Real$ is star shaped with
respect to $x \in S$ then $S$ must necessarily be star shaped
with respect to $x$. Alternatively, we can use the following
definition.
\begin{dfn}\label{a:star_shaped_func}
A function $f:S\rightarrow\Real$ (on a star shaped  set $S$
with respect to $x$) is star shaped with respect to $x$ iff the
following holds for any $y \in S$ and every $\gamma \in [0,1]$:
\begin{equation}\label{eq:new_Jensen}
f(\gamma x + (1-\gamma) y) {\leq} \gamma f(x) + (1-\gamma) f(y).
\end{equation}
\end{dfn}
Expression (\ref{eq:new_Jensen}) is a form of Jensen's inequality with one ``fixed end''.
It is obvious that a function on $S$ is convex iff it is star
shaped with respect to every point $x\in S$.

\begin{dfn}\label{a:star_shaped_envelop_func} Let
$f:S\rightarrow\Real$ be a bounded from above and below
function ($A < f(x) < B$). The supremum of  star-shaped (with
respect to $x$) minorants of $f(\cdot)$ is the star shaped envelope of
$f(\cdot)$ with respect to $x$, ${\rm star}_x (f)$.
\end{dfn}
Note that Definition \ref{a:star_shaped_envelop_func} can be viewed as a definition of the star-shaped envelope of a function ''from  above" (compare with (\ref{a:star_shaped_envelop_set_above})).

Let  $S \subset V$ be a star shaped set with respect to $x \in
S$,   and ${\rm conv} (f)$ be the convex envelop of $f(\cdot)$. The following properties hold for the star shaped sets and
functions:
 ${\rm conv} (f)(y) \leq {\rm star}_x (f)(y) \leq f(y)$.
The first property follows immediately from the definition, and
the second property is a consequence of the following fact:
\begin{equation}\label{a:conv-star}
{\rm epi}(f)\subset{\rm epi}({\rm star}_x (f))\subset{\rm epi}({\rm conv} (f)).
\end{equation}
Let us follow the definition of a star envelope of a set from
below, (\ref{defstarenvset}), and produce the definition of a
star envelope of a function ''from below". Consider the
one-dimensional case: $V=\Real$. Let $S=[a,b]$, $x \in S$; for
any $z\in S$ ($z \neq x$) we define
\begin{equation}\label{a:minorant}
\phi_{z,x}(y)=\left\{
 \begin{aligned}
 & \min\left\{f(y), \frac{y-x}{z-x}f(x) + \frac{z-y}{z-x}f(z)\right\} \; \mbox{if} \; 0< \frac{y-x}{z-x} <1 ; \\
 & f(y) \; \mbox{else}.
 \end{aligned} \right.
\end{equation}
As follows from Jensen's inequality with one fixed end, (\ref{eq:new_Jensen}), the
following holds: $${\rm star}_x(f)(y)=\inf_{z \in S} \{
\phi_{z,x}(y)\}.$$
One can also see that the definitions of star-shaped envelopes of a function ``from above'' and ``from below"  are equivalent.

The following proposition is obvious:

\begin{prop}\label{prop:star_properties} Let $S=[a,b]$ be a closed interval in $\Real$, $p$ be an element from $S$, and $f:S\rightarrow\Real$ be a continuous function. Consider the functions $\phi_{z,p}(\cdot)$, $z\in S$, defined as in (\ref{a:minorant}). Then

1) the functions $\phi_{z,p}(\cdot)$, $p\in S$ are equicontinuous
with $f(\cdot)$;

2) the set of functions $\{ \phi_{z,p}(\cdot)\}$, $p\in S$ is
compact;

3)  the function ${\rm star}_p(f)$, $p\in S$ exists and is a
continuous function.

4) if a continuous function on $S$ achieves its minimum at a
single point $x$, then its star shaped envelope with respect to
$x$ has the same property.

5) the star shaped envelope (w.r.t. a point $p\in S$) of a
monotone function  is monotone.
\end{prop}

The first property follows immediately from (\ref{a:minorant}).
Properties 2 and 3 follow from the Arzela--Ascoli theorem. Indeed,
according to this, equicontinuity and uniform
boundedness of $\phi_{z,p}(\cdot)$, $p\in S$, imply that
$\{\phi_{z,p}(\cdot)\}$ is relatively compact. Compactness the
follows from the fact that the set of functions
$\{\phi_{z,p}(\cdot)\}$ is closed.  To demonstrate existence and
continuity of ${\rm star}_p(f)$, $p\in S$ consider a
sequence $\{g_i\}_{i=1}^\infty$ of grids $g_i=\{a,a+(b-a)/i,
a+(b-a)k/i,\dots,b\}$, $1\leq k< i$ on $S$, and define
\[
f_i(y)=\inf_{r\in g_{2i}}\{\phi_{r,p}(y)\}
\]
The sequence of functions $\{f_i(\cdot)\}_{i=1}^{\infty}$ is
equicontinuous and $f_1(\cdot)\geq f_2(\cdot)\geq\dots \geq f_n(\cdot)\geq \dots$
pointwise. This means that the sequence converges uniformly,
and that $\lim_{i\rightarrow\infty}f_i(y)=\inf_{z \in S} \{
\phi_{z,p}(y)\}={\rm star}_p(f)(y)$. Continuity of the limiting
function, ${\rm star}_p(f)$, follows from equicontinuity of
the family $\{f_i(\cdot)\}_{i=1}^\infty$.  Property 4 is easily
verifiable by the contradiction argument. Property 5 is the
consequence of that the functions $f_i(\cdot)$ are monotone (by
construction) if the function $f(\cdot)$ is monotone. $\square$

Let $V=R^n$, and $S\subset V$ be a compact and star shaped with
respect to $x$. For every $z\in S$, and $\gamma \in [0,1]$ we
define $\psi_z(\gamma)=\min\{f((1-\gamma)x + \gamma z),
(1-\gamma)f(x) + \gamma f(z) \}$. Then

$${\rm star}_x(f)(y)=\inf_{(1-\gamma)x + \gamma z = y, \, \gamma \in [0,1], \,  z \in S} \{ \psi_z(y)\}.$$

The properties of ${\rm star}_x(f)(y)$ depend on the
properties of mapping
\begin{equation}\label{auxmapp}
y \mapsto \{(\gamma,z) \mid   (1-\gamma)x +
\gamma z = y, \, \gamma \in [0,1], \, z \in V\}.
\end{equation}
If this mapping is continuous (in the Hausdorff metrics in the
space of compact sets), then the star shaped envelope of every
continuous in $S$ function with respect to $x$ is also
continuous. Mapping (\ref{auxmapp}) is continuous in $S$ iff
the Minkovski functional
$$p(y)=\inf_{a>0}\{a \mid y-x \in a(S-x)\}.$$
is continuous in $S$.

\end{document}